\documentclass[12pt]{amsart}
\usepackage{amsmath}
\usepackage{amssymb}
\usepackage{graphicx}

\newfont{\bigbf}{cmbx10 scaled\magstep1}
\normalbaselines
\def\C{\mathbb C}
\def\R{\mathbb R}
\def\H{\mathbb H}
\def\Z{\mathbb Z}

\def\E{\mathbb E}
\def\ka{K\"ahler}
\def\hc{\H_{\C}^2}

\def\uo{\operatorname{U}(1,1)}

\theoremstyle{definition}

\newtheorem{dfn}{Definition}[section]
\newtheorem{defn}[dfn]{Definition}

\newtheorem{rem}[dfn]{Remark}

\theoremstyle{plain}

\newtheorem{thm}[dfn]{Theorem}

\newtheorem{lem}[dfn]{Lemma}

\newtheorem{prop}[dfn]{Proposition}

\newtheorem{cor}[dfn]{Corollary}

\newtheorem{conj}[dfn]{Conjecture}

\def\diam{\operatorname{diam}}
\def\ind{\operatorname{ind}}

\title{Intersection pairing  on
 hyperbolic 4-manifolds}
\author{Michael Kapovich}
\date{\today}

\begin{document}

\maketitle

This paper was written as an MSRI preprint in 1992. For a long time I was planning to revise and clean it up. Instead, I am just posting the paper ``as is.''

\section{Introduction }

\bigskip

1.1 We begin with the following elementary question:
How one can estimate the algebraic intersection number 
$< \alpha , \beta >$ between two 
closed curves $\alpha$, $\beta$ on a complete oriented Riemannian surface $M$?
 
Suppose that the metric on $M$ is hyperbolic (has curvature $(-1)$).	
Then the answer can be given in the terms of lengths of $\alpha, \beta$.
For example:
$$
|< \alpha , \beta >|  \le K(l( \alpha ) , l( \beta ))= 
2e^{l(\beta )}(\pi /2  +   l(\alpha )) 
+ \frac{\pi/2 \ + l(\alpha )}{l(\beta ) /2  - \log (\exp (l(\beta ) /2) -1)}
 $$
In particular, if $|<\alpha , \beta >| \ge 1$ then 
$$ 
1 \ \le \ \sinh (l(\beta ) /2) \ \sinh (l(\alpha ) /2)
$$

Note that the right side depends on the metric 
introduced on $M$  while the left side is purely topological.
 
\medskip 
1.2. Now let $M$ be an arbitrary complete oriented hyperbolic 4-manifold,
 $\sigma _j : \Sigma _j \to M \ (j= 1,2)$ be two cycles in $Z_2(M, \Z)$,
 where $\Sigma _j$ are closed oriented connected 
surfaces. The main aim of this paper
 is to estimate from above the absolute value of the intersection pairing
$| <[\sigma _1  ] , [ \sigma _2 ] > |$ so that the estimate depends
just on the Euler characteristics of  $\Sigma _j$. This will be done
under certain  condition on  $\sigma _j$.

\begin{defn}
The maps $\sigma _j$ are { \em incompressible} if the induced 
homomorphisms of fundamental
groups $\sigma _{j*}$ are injective.
\end{defn}

\begin{thm}\label{t1}
%THEOREM 1. 
There exists a function $C(.,.)$ such that for any complete
hyperbolic 4-manifold $M$ and for any 
 classes $[\sigma _1  ] , [ \sigma _2 ]$
 in $H_2(M, \Z )$ which have incompressible 
representatives, 
we have:
$$
 | <[\sigma _1  ] , [ \sigma _2 ] > | \le C(| \chi (\Sigma _1  ) | , |\chi (
 \Sigma _2 )|) \ .$$
\end{thm}

\begin{thm}\label{t1'}
%THEOREM 1'. 
The conclusion of Theorem \ref{t1} holds even if instead of 
the {\em algebraic intersection number} $< [\sigma _1  ] , [ \sigma _2 ] >$
 we take the {\em geometric intersection number} 

$i( \sigma _1 , \sigma _2 ) := \inf \{$ { \em number of points of intersection
between} $\varphi _1 ( \Sigma _1 ) , \varphi _2 ( \Sigma _2 ) :
\varphi _j $ {\sl is homotopic to} $\sigma _j , j= 1, 2 \}$  
\end{thm}

Thus, Theorem \ref{t1} shows how the intersection pairing together with simplicial
norm on the second homology group provides an obstruction for existence of 
complete hyperbolic structure on
 4-manifolds. The intersection pairing itself is not  too
interesting invariant for hyperbolic 4-manifolds. If $M^4$ is closed
then just two symmetric bilinear forms of the given rank can occur.
On another hand, any  symmetric bilinear form can be realized as
 the intersection form of compact convex  hyperbolic  4-manifold. 
There are no homological obstructions for hyperbolicity in the dimensions
2 and 3  because the length of curve (unlike the area of surface) 
is not  a topological
invariant.

1.3. Consider the special case: suppose $M= M(e,g)$ is 
homeomorphic to a $\R ^2$- bundle over a closed orientable surface $F$ 
of genus $g > 1$. Such bundles are characterized by their "Euler
 numbers". In this case it is just the "selfintersection" number
 $e= <[ F ], [F] >$, where we identified $F$ with the zero section of 
the bundle. 

 \begin{cor}
  \label{c1} 
  The condition $|e| \le C(2g-2, 2g-2)$ 
is necessary for existence of complete hyperbolic structures on $M$.
\end{cor}

%REMARK 1. 
%\begin{rem}\label{r1}
%It follows from the calculations in Section 3, as 
%$C(2g- 2, 2g- 2)$ one can actually take 
%$$
 % \exp (8000 g/ \mu )$$
%where $\mu$ is the Margulis constant for the 4-dimensional hyperbolic
%space.
%\end{rem}

Denote by ${\mathcal S}(e, g)$ an orientable 3-manifold 
which is a circle bundle over the closed oriented surface of genus $g$
such that the Euler number of the fibration is $e$.

\begin{cor}\label{c2}
%COROLLARY 2. 
 The condition $|e| \le C(2g-2, 2g-2)$ 
is necessary for existence on ${\mathcal S}(e, g)$ flat conformal structures   
with nonsurjective development maps.
\end{cor}

%For another example of application of Theorem 1' consider a pair of
%Kleinian groups $G_1 , G_2 \subset SO(4, 1)$ such that $G_j$ are isomorphic
%to $\pi _1 (S_{g_j} )$ (where $S_g$ is a closed orientable surface of genus $g$);
% and the limit sets $L_j = \Lambda (G_j )$ are topological unknot in
%$\S ^3$ ($j = 1, 2$). Suppose that the group $G$ generated by $G_1 , G_2$ is
%discrete and torsion-free. If we do not control the genera $g_j$ then $L_1 , L_2$
%can form an arbitrary link in $\S ^3$ (see for example [GLT]). However,
%Theorem 1' implies that given a pair $(g_1 , g_2 )$ there exist not more 
%than $k(g_1 , g_2 )$ links $(L_1 , L_2 )$ such that the group $G$ is 
%discrete and torsion-free, where  
%$k(g_1 , g_2 ) \le \exp (8000 (g_1 + g_2)/ \mu )$.

1.4. The examples of hyperbolic manifolds $M= M(e,g)$ with $e=0$ are 
 easy to construct and they  exist for arbitrary genus $g > 0$. If
the surface $F$ of genus $g$ is uniformized by the group $\Gamma \subset 
PSL(2, \C )$
then the extension of  $\Gamma $ in $\H ^4$ is the holonomy group for
the complete hyperbolic structure on $M(0, g)$. As the
application of the Bieberbach theorem we obtain that for $e > 0$
 manifolds $ M(e,1)$ can't
have complete hyperbolic structure.
 The case $g > 1 \ , \ e > 0$ is less trivial, first  examples of
complete hyperbolic structures were obtained in \cite{GLT,Ku1,Ka1} (with the estimate $(2g-2)/3$ due to Kuiper):

\begin{thm}\label{t2}
%THEOREM 2. 
The manifold $M= M(e,g)$ admits a complete 
hyperbolic structure assuming that 

%(Kuiper)
$$
0 < e \le (2g-2)/3  $$ 
%(Kapovich) 
%$$
%0 < e \le (2g-2)/22 \ .$$
\end{thm}

After \cite{GLT} M.Anderson, \cite{A}, proved

%THEOREM 3. 
\begin{thm}\label{t3}
Let $E \to B$ be an arbitrary vector bundle 
with the compact base $B$ of negative sectional curvature. 
Then $E$ admits a  Riemannian metric
with strictly 
negatively pinched sectional curvature:
$$
0 > a_E > K_E > -1$$
for some constant $a_E$ depending on the bundle.
\end{thm}

%CONJECTURE 1
\begin{conj}\label{con1}
(See \cite{GLT}.) The (Milnor-Wood) inequality 
$$
0 \le |e| \le (2g-2)$$
is the necessary condition for existence of complete hyperbolic
 structure on the manifold $M= M(e,g)$.
\end{conj}

%REMARK 2. 
\begin{rem}\label{r2}
It is important here that $M$ is fibered.
 N.Kuiper, \cite{Ku2}, constructed a sequence of complete hyperbolic manifolds $M_g^4$
 which are homotopy equivalent to closed surfaces $F_g$ of genus $g$ 
such that:
$$
\lim_{g \to\infty} 
<[F_g], [F_g] >/(2g-2) = 2/ \sqrt {3} >  1 $$
\end{rem}

%REMARK 3
\begin{rem}\label{r3} (N.Kuiper, \cite{Ku3}). 
Suppose that $M$ is hyperbolic and
$\Sigma \subset M$ is an {\em imbedded} minimal surface of genus $g$.
 Then the Milnor-Wood
 inequality
$$
 |<[\Sigma ], [\Sigma ]>| \le (2g-2)$$
 holds.
\end{rem}

1.5. For another example of application of  Theorem \ref{t1'}  consider a pair of
Kleinian groups $G_1 , G_2 \subset SO(4, 1)$ such that $G_j$ are isomorphic
to $\pi _1 (S_{g_j} )$ (where $S_g$ is a closed orientable surface of genus $g$);
 and the limit sets $L_j = \Lambda (G_j )$ are topological unknots in
$S^3$ ($j = 1, 2$). Suppose that the group $G$ generated by $G_1 , G_2$ is
discrete and torsion-free. If we do not control the genera $g_j$ then $L_1 , L_2$
can form an arbitrary link in $S^3$ (see for example [GLT]). However,
 Theorem \ref{t1'}  implies that given a pair $(g_1 , g_2 )$ there exist not more 
than $k(g_1 , g_2 )$ links $(L_1 , L_2 )$ such that the group $G$ is 
discrete and torsion-free, where  
$k(g_1 , g_2 ) \le \exp (8000 (g_1 + g_2)/ \mu )$.

1.6. Probably it's possible to prove  Corollary \ref{c1}  for  convex compact 4-manifolds 
 by comparing two $\eta $-invariants
for flat conformal manifolds, \cite{Ka2, Ka3}. A more realistic idea was suggested
to author by M.~Gromov who proposed to compactify the moduli space of all hyperbolic
structures on the given fiber bundle. Formally speaking this idea does not
work, since arbitrary large number of self-intersections of zero section
can be pinched to point in limit.  However, what we are using in this paper
are ``pre-limit" considerations based on Mumford's compactness theorem and
existence of Margulis constant. 

The idea of the proof is quite simple. Suppose
that we realized the classes $[\sigma _1], [\sigma _2]$ by piecewise-geodesic surfaces 
$\Sigma _1, \Sigma _2 $ in $M$ such that the number and diameter of simplices in $\Sigma _1 \cup \Sigma_2 $ are
 bounded from above. 
 
 Then the existence of the universal Margulis constant
and the fact that two geodesic planes  intersect transversally in not more
than one point immediately imply the assertion of  Theorem \ref{t1} . Certainly it's impossible in
general to estimate from above the diameter of simplices since the diameters
of the surfaces $\Sigma _1, \Sigma _2 $ can be unbounded. However the ``long"
pieces of 
$\Sigma _1, \Sigma _2 $ are contained in the ``thin" part of the manifold $M$
which have very simple topological structure. Then the detailed analysis of
the situation in the ``thin" part of $M$ (section 2) and correct choice of the
surfaces representing $[\Sigma _j]$ (sections 3, 4) yield the claimed result.

1.7. Thus,  Theorem \ref{t1}  is a ``0-th order approximation" to  Conjecture 1.
The simplest examples of  negatively pinched closed manifolds of 
dimension $4$ 
which do not admit hyperbolic structure are given by complex-hyperbolic 
manifolds. More sophisticated examples were constructed by Mostow and Siu, 
\cite{MS};  there are no (real) hyperbolic structures on any compact 
K\"ahler manifold of (real) dimension $>2$
 due to theorem of J.Carlson and D.Toledo, \cite{CT}. Another series of
examples was constructed by Gromov and 
Thurston in  \cite{GT}.  Theorem \ref{t1}  combined with the theorem of Anderson presents 
the first examples of negatively curved  open manifolds of dimension 4
that do not have complete hyperbolic structure. It's interesting to remark
that these manifolds have finite-sheeted branched  coverings which are 
hyperbolic (as well as some examples of Gromov and Thurston).

Certainly there is a gap between Theorems 1 and 3 and Conjecture 1. 

%CONJECTURE 2. 
\begin{conj}\label{con2}
There exists a function $D(.,.,.)$ such that for any
Riemannian 4-manifold $M$ whose curvature is pinched as $0 > a \ge K_M \ge -1$
 and for any classes $[Q]$, $[P]$ in $H_2(M, \Z )$ we have:
$$ 
\vert <[Q], [P] > \vert \ \le D(\| [Q] \| , \| [P] \|, a),$$
where $\|W \| = \min \{| \chi(W) |; w: W \to M $ is a
 surface representing the class $[W]
\}$.
\end{conj}

%REMARK 4 
\begin{rem}\label{r4} 
(W.Goldman). Conjecture 2 is not true for {\em orbifolds}. The
example is given by complex-hyperbolic orbifolds covered by nontrivial
 $\R ^2$-bundle over surface. Namely, let $\Gamma \subset SU(1,1)
 \subset SU(2,1)$ be a cocompact torsion-free lattice in $SU(1,1)$.
Then the manifold $M(\Gamma ) = \H_{\C}^2 / \Gamma$ is complex- hyperbolic
and it admits an isometric $U(1)$-action. On other hand, $M(\Gamma )$ is
diffeomorphic to the total space of a nontrivial $\R ^2$-bundle over 
$S= \H^2_{\R} /\Gamma$ (2-sheeted ramified covering over the tangent bundle
of $S$). The group $U(1)$ has cyclic subgroups
 $\mathbb{Z}_n$ of arbitrarily large order $n$. Then the sequence of orbifolds 
$M(\Gamma )/\mathbb{Z}_n$ has the desired properties.
\end{rem}

As the particular case we have:

%CONJECTURE 2'.  
\begin{conj}\label{con2'}
In  Theorem \ref{t1}  one can drop the condition
for the cycles to be incompressible.
\end{conj}

1.8. If true, Conjecture \ref{con2'} would have several applications for flat
conformal structures on 3-manifolds described below.

%Denote by ${\mathcal S}(e, g)$ an orientable 3-manifold 
%which is a circle bundle over the closed oriented surface of genus $g$
%such that the Euler number of the fibration is $e$.

%COROLLARY 2.  The condition $|e| \le C(2g-2, 2g-2)$ 
%is necessary for existence on ${\mathcal S}(e, g)$ flat conformal structures   
%with nonsurjective development maps.

%Let  $N$  be any closed Haken 3-manifold which is obtained by gluing of
%geometric manifolds. Flat conformal
%structure $K$ on $N$ is called to be {\em extendable}
%if $(N, K)$ can be realized as an ideal boundary component of a complete
%hyperbolic 4-manifold $M^4$. The extendable structure is called to be 
%{\sl nice} if the inclusion $N \to M^4$ induces a monomorphism of the 
%fundamental group.
%it has a  nonsurjective development map
%and the kernel of the holonomy representation
%$$
%\rho : \pi _1(N) \to SO(4, 1)$$
% is the normal closure of homotopy classes of regular fibers
%of Seifert components. (So $Ker(\rho )$ is as small as the discreteness of
%$\rho ( \pi _1(N))$ allows.)

Let $N_1$ and $N_2$ be  compact oriented geometric 3-manifolds such
that:

(1) interiors of $N_i$ have no Euclidean structures;
 
(2) $\partial N_i = T_i$ are incompressible tori $ (i= 1, 2)$.

Denote by $N_1 \cup _f N_2$ the manifold obtained by gluing of $N_1 , N_2$
via the homeomorphism $f : T_1 \to T_2$.

%COROLLARY 3. 
\begin{cor}\label{c3}
(i) If $N_i$ are both hyperbolic then not more than
finitely many manifolds $N_1 \cup _f N_2$ can be realized as incompressible
ideal boundary components of a complete
hyperbolic 4-manifolds. (This is a consequence of Morgan's compactness
theorem, \cite{Mor}).
\end{cor}
%admit a nice extendable 
%flat conformal structure.

If Conjecture \ref{con2'} holds, then we have the conclusions (ii) and (iii)

 (ii) If both manifolds are Seifert then not more than
finitely many manifolds $N_1 \cup _f N_2$ can be realized as
ideal boundary components of  complete
hyperbolic 4-manifolds.
%admit an extendable 
%flat conformal structure.

(iii) If (say) $N_1$ is hyperbolic, $N_2$ is Seifert then  infinitely
many manifolds $N_1 \cup _f N_2$ can be realized as
 ideal boundary components of  complete
hyperbolic 4-manifolds (cf. \cite{Ka1, Ka4}).
%admit an extendable 
%flat conformal structure.

However, if we fix the image of the regular fiber of $N_2$ then again 
there are only finitely many manifolds $N_1 \cup _f N_2$ which can be 
ideal boundary components of  complete hyperbolic 4-manifolds.

%with an extendable 
%flat conformal structure.
% then there
%is not more than finitely many isotopy classes of homeomorphisms 
%$f : T_1 \to T_2$ such that:

%the closed oriented manifold $N= N_1 \cup _f N_2$ admits a minimal
%uniformization with torsion-free holonomy group.

%(iii) If (say) $N_1$ is hyperbolic, $N_2$ is Seifert then there are infinitely
%many isotopy classes of $f$ with the property of minimal uniformization.
%However,
%if we fix the image of the regular fiber of $N_2$ then again there are only
%finitely many possibilities for the minimal uniformization.

On other hand it follows from \cite{Ka1, Ka4} that unless the canonical
decomposition of the Haken manifold $N$ includes gluings of the type  (i) 
the manifold $N$ always has a finite-sheeted covering $N_0$
which is an ideal boundary component of a complete hyperbolic 4-manifold.
%admits an extendable flat conformal structure .

\bigskip
 
1.9. We split the proof of  Theorem \ref{t1}     in two cases.

Case 1 (Section 3). We shall suppose that the cycles $\sigma _j$  satisfy 
 certain  condition of maximality.

{\bf Condition MAX.}  If for some $h \in \pi _1(M)$ we have
 $$1 \ne g  \in h^{-1} \sigma _{1*}( \pi _1 (\Sigma _1) )h \cap \sigma _{2*}
( \pi _1 (\Sigma _2))  $$
then the maximal  virtually abelian subgroups of
$$
\pi _1(M),  h^{-1} \sigma _{1*}( \pi _1 (\Sigma _1) )h \ , \sigma _{2*}
( \pi _1 (\Sigma _2)) $$
 containing $g$ are equal.
(For example, in Theorem 2 the condition "MAX" is satisfied.)

Case 2 (Section 4). This is the general case.

In the Case 1 the analysis of behavior of the surfaces
 $\sigma _j (\Sigma _j )$ in the ``thin
part" of the manifold $M$ is simpler, which is why we decided to single
out this case.

In the section 2 we discuss the geometry of components of
 the ``thin" parts of hyperbolic
4-manifolds (``Margulis tubes and cusps"). The following is the reason of 
difference between dimension 3 and 4. Let $\langle g\rangle$ be an infinite
 cyclic group of (orientation preserving) isometries of $\H ^n$.
Consider the set of points ${\mathcal K} (\langle g\rangle , \mu ) = \{ x \in \H ^n :
d(x, g^k (x)) \le \mu $ for some $\mu \ne 0 \}$; define $q(x) := $ minimal
$k > 0$ such that $d(x , g^k (x)) \le \mu $. Then, for $n \le 3$ the function
$q$ is constant on  ${\mathcal K} (\langle g\rangle, \mu )$. However it's not longer true
for $n \ge 4$. It is well-known that the most dramatical situation is in the 
case of parabolic 
$g$, when $q(x)$ can have infinitely many different values. If $g$ is
loxodromic, then the image of $q$ is still finite, but it depends on
the 
%"rotational part" of the 
element $g$. Something similar occur even
 for $n \le 3$ if $g$ does not
preserve the orientation; however, in this case $q$ can have not more
than 2 different values. 

1.10. There are several other results and conjectures that seems to be
similar to
 Theorem \ref{t1}  and Conjecture 1. 

1.10.1. If $C$ is a smooth curve of genus $g$
in a complex surface $X$ and $K$ is a canonical class of $X$ then
$$
2g - 2 = < C, C > + K \cdot C$$
(see \cite{F}).

R.Kirby conjectured, \cite{Ki}, that for any smooth embedded 2-manifold $\Sigma$
in the same homology class as $C$ then genus of $\Sigma$ is $\ge g$.

Conjecture of J.Morgan is that for any smooth oriented 4-manifold for which
Donaldson's polynomials are defined and non-zero, and any smoothly 
embedded oriented surface $\Sigma \subset M$ with positive self-intersection
one have the inequality

$$
2g - 2 \ge \ \ <\Sigma , \Sigma >$$

I learned this information from the paper of P.Kronheimer [Kr] where
the reader can find further information on this subject.

1.10.2.  Suppose that $M$ is a complex hyperbolic surface and
$f: S_g\longrightarrow M$ is a homotopy-equivalence.
Let $\omega_M$ denote the \ka\ form on $M$.
Then Domingo Toledo has proved in \cite{To1} that the number
$$
c = \frac{1}{2\pi} \int_\Sigma f^*\omega_M $$
is an integer independent of $f$ which satisfies
$$
2-2g \le c \le 2g-2  $$
Furthermore Toledo, \cite{To2}, proved that $M$ is a quotient by a cocompact lattice
in $\uo$ if and only if $\vert c\vert = 2g-2$.

In \cite{GKL} we proved that, subject to Toledo's necessary
conditions, every value of $c$ is realized by a complex hyperbolic
surface $N( c, g)$ homeomorphic  to $M(e=e(g, c), g)$ (see 1.3). 
In all these examples the
self-intersection number $e=e(g, c)$ of the generator of $H_2 (N(c, g) , \Z )$
varies in the closed interval $[1 -g, 2(1-g) ]$. So, in particular, some
of these manifolds are homeomorphic, but the actions of their fundamental
groups on $\hc$ can't be deformed one to another inside the group
Isom($\hc$). On other hand, W.Goldman showed that
 in the examples \cite{Ka1} (see  Theorem \ref{t1} ) all representations
of $\pi _1 (S_g )$ in $SO(4, 1)$ are in the component of the 
trivial representation in Hom($\pi _1 (S_g ), SO(4, 1)$).

1.10.3. See also the paper of N.Mok, \cite{Mo}.

\bigskip

1.11. {\bf Acknowledgements.} I am deeply grateful to Misha Gromov and
Nicolaas Kuiper for reviving my interest to the subject of the current paper
(Conjecture 1), helpful advises and discussions. 
This work was supported by NSF grants numbers 8505550 and 8902619 
administered through the University
of Maryland at College Park and MSRI which the author gratefully acknowledges.

\section{Geometry of Margulis cones}

\bigskip
Many results of this section are well known in some form, 
some are known as folklore.

2.1. {\bf Notations.} By $d(x , Y)= \inf \{ d(x, y) : y \in Y \}$ we shall denote the
(low) hyperbolic 
distance between the 
1-point set $\{ x \} $ and the set $Y \subset \H ^4$. The ball with the center
at $x$ and
radius $r$ will be denoted by $B(x, r)$.  If $h \in Isom(\H ^4)$ then
 $l(h) = \inf \{ d(h(x), x) : x \in
 \H ^4 \} $
is the ``length" of $h$. For each pair of points $a, b \in \H ^4$ we shall denote
by $[a, b]$ the geodesic segment connecting them. We denote by $[a, b , c]$
 the totally geodesic triangle 
with the vertices $a, b , c \in \H ^4$. If $G \subset Isom(\H^4)$, $x \in
\H ^4$ then
$Ir_G(x) = d(x, Gx)/2$ is the injectivity radius $InjRad([x])$
at the projection $[x]$ of $x$ in
$\H ^4/G$.
We shall assume that all groups below are torsion-free.

 If $h$ is a loxodromic or parabolic transformation in $\H ^4$ then we 
denote by $\Pi $ the canonical fibration of $\H ^4$ by totally geodesic 
hyperplanes orthogonal to either axis of $h$ (if $h$ is loxodromic) or 
its 1-dimensional invariant horocycle (if $h$ is parabolic). 
%If $h \in H$-
% virtually abelian discrete group then 
The projection of $\Pi $ to $\H ^4 / \langle h\rangle$ will be  called {\em canonical 
foliation}
associated with $\langle h\rangle$.

 For the  virtually abelian group
$H \subset $Isom$(\H ^4)$ define 
$$ {\mathcal K} (H, \mu ) = \{z \in \H ^ 4 : Ir_H(z) \le \mu \} \eqno(1)$$
to be the {\em Margulis cone}. The projection 
$$
{\mathcal T}(H, \mu )= {\mathcal K}(H, \mu ) /H \eqno(2)$$
 is the {\em Margulis tube} in $\H ^4 /H$ (we assume that Margulis tubes can
be noncompact; in particular cusps are also considered as  Margulis tubes).

Suppose that the hyperbolic 4-space is realized as "upper half-space",
and the loxodromic element $h$ is 
%has one fixed point at infinityof $\R ^4$.
a Euclidean similarity. The main problem concerned with
Margulis tubes in the dimension 4 (and higher) is that even for the cyclic
loxodromic group
$\langle h\rangle = H \subset $Isom$(\H ^4)$ the boundary of
the Margulis cone ${\mathcal K}(H, \mu )$ is very far from been the "round" cone 
(like in
dimensions 2 and 3), but
 rather
looks as a cone over an ellipsoid, where the ratio of largest and smallest
 axes
can be arbitrary large.

%Many statements below are well-known and we present them only for
%completeness.

\bigskip
2.2.	

%LEMMA 1.

\begin{lem}\label{l1}
 Let $x \in \H ^4$ be a point such that: for some subgroup
$H \subset $Isom$(\H ^4)$ we have: $Ir_H(x) \ge \nu$ for some positive 
$\nu$. Then the ball $B(x, r)$ contains not more than 
$$
(  \exp^3(r + \nu ))/\nu^3 \eqno(3)$$
points from the orbit $Hx$.
\end{lem}

\begin{lem}\label{l2}
	%LEMMA 2. 
	Under the assumptions of Lemma \ref{l1}, the number of elements
 $h \in H \subset $Isom$(\H ^4)$ such that the intersection
$h(B(x,r)) \cap B(x,r)$ is non-empty is at most 
$$
(  \exp^3(2r + \nu ))/\nu^3 \eqno(4)$$
\end{lem}

\begin{lem}\label{l3}
	%LEMMA 3. 
	Let $H$ be a discrete subgroup of Isom$(\H ^4)$ and 
$\nu /2 = Ir_H(x)$. Suppose that $d(x, y) < r$. Then 
 $Ir_H(y) > C_1(r, \nu )/2$. Here 
$$ 
C_1(r, \nu )= 2r/n(r, \nu), \ n(r, \nu ) = \left[\frac{\exp (18r + 2 \nu)}{\nu ^3} \right] + 1 \eqno(5)$$

\end{lem}

\proof 
Let $h$ be arbitrary nontrivial element of the group $H$.
 Let $n_0$ be such that $d(x, \ h^{n_0}(x)) \ge 3r$. Then $d(y, \ 
h^{n_0}(y)) 
\ge r$ and $d(y, \ h(y)) \ge r/n_0$. So, our aim is to estimate this $n_0$
from above. Notice that for $n = n(r, \nu )$ among the elements
$$
\{ 1, \ h, \ ... \ , h^n \}$$
there is $h^k$  such that $d(x, \ h^k(x)) \ge 3r$ (by Lemma 1). Then we can
take 
$n_0 \le n$ and $d(y, h(y)) \ge r/n $ for every $h \in H - \{ 1 \} $.

Lemma follows. \qed

\begin{figure}[htbp] %  figure placement: here, top, bottom, or page
   \centering
   \includegraphics[width=4in]{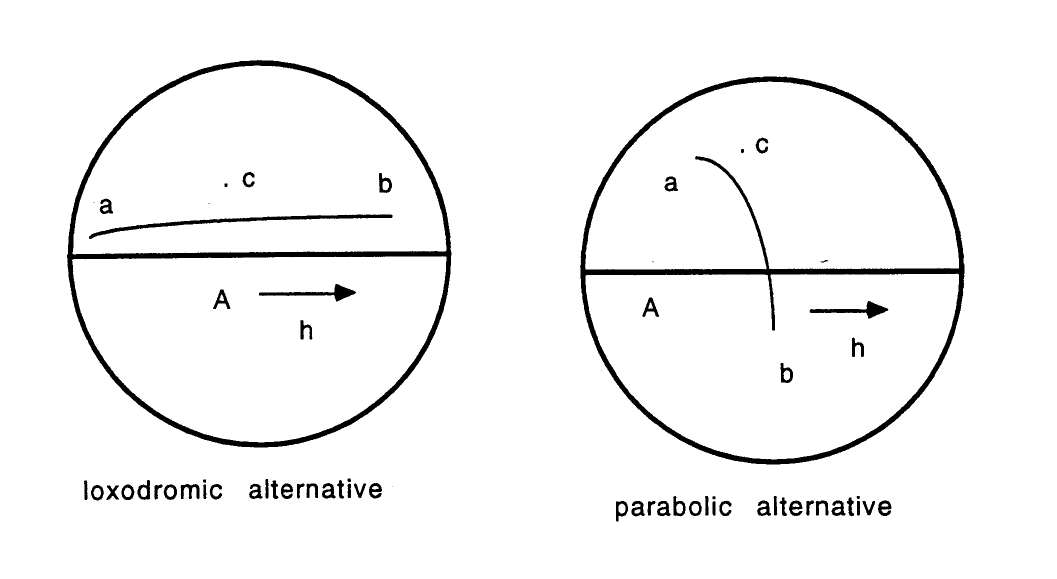} 
   \caption{}
   \label{F1}
\end{figure}

2.3. In Lemmas \ref{l4} and \ref{l5} below we shall prove the {\em Key Property of Index}  
(Corollary \ref{c5}) which will be crucial to our paper: 

{\bf Key property of index.} Let $h$ be a parabolic or loxodromic
isometry of $\H^4$,  $x \in \H^4$. The index $ind_h(x)$ of $x$ 
 is just $d(x, hx)$.
 Suppose that $x$ is such that: $Ir_{\langle h\rangle} (x) > \mu $ and 
there is a geodesic segment  $L = [a, b]$ such that $d(x, L)
 < R $ for some  $R$  and
$$
\max \{\ind_h(x), \ind_h(a), \ind_h(b) \}  < R.$$ 
Then we shall prove that 
 either  
$$
\min \{ d(x,a), d(x,b)\} < C_+(R, \mu )$$
 (parabolic alternative) or 
$$
l(h) > C_-(R, \mu ) > 0 $$
 (hyperbolic alternative). See Figure \ref{F1},
where $A= A_h$ is the axis of the loxodromic transformation $h$.

\begin{rem}
%REMARK. 
The Key Property of Index is valid also if instead of $\H ^n$ we shall take
arbitrary simply connected complete space $X$ of sectional curvature
$K_X < -a^2 < 0$.
\end{rem}

We shall prove the Key Property of Index only for the loxodromic $h$,
 the parabolic case easily follows.

\bigskip

Below we assume that the hyperbolic 4-space $\H ^4$ is realized as the 
upper half
space $\R ^4 _+$  ; \ $|X - Y|$ denotes the Euclidean distance between
points $X , Y$; $ \partial_{\infty } \H ^4 = \overline{\R ^3} = \R ^3 \cup
\infty$.

\begin{figure}[htbp] %  figure placement: here, top, bottom, or page
   \centering
   \includegraphics[width=4in]{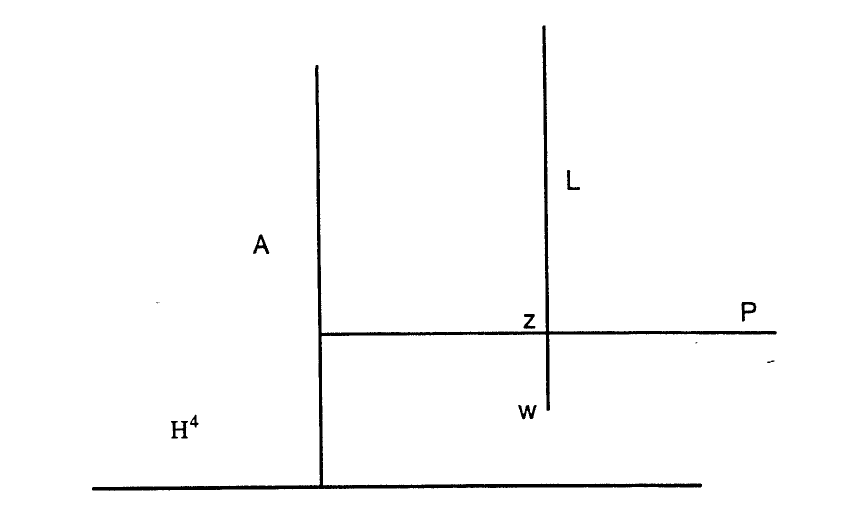} 
   \caption{}
   \label{F2}
\end{figure}

2.4.  Suppose that $g = \theta \circ \lambda $ is  a similarity in
 $\E ^4$ 
preserving $\H ^4$,
$g(0)=0, \ A$ is  the axis of $g$, $d(z, A) > 2$. Here $\lambda$ is
the dilatation $\lambda : x \to \lambda x$ and $\theta$ is the rotation on
the angle $\theta$. Let $L$ be the geodesic
containing the points $\infty , z$; let $w \in L$ be a point such that
$z$ lies between $w$ and $\infty$. (See Figure \ref{F2}). 

\begin{lem}\label{l4}
%LEMMA 4. 
Suppose that under conditions above: 
$$
\nu \le d(g(z), z) \le R ;\  d(g(w), w) \le R.$$
Then $d(z, w) \le R + \frac{1 }{ \nu}$.
\end{lem}

\proof Step 1.  Set $g(z)= z' \ , g(w) =w'$. Denote by $\alpha (u)$ the angle
 between the horosphere $P$ with center at
$\infty$ containing  the point $z$ and the Euclidean line passing through
the points $z, u$, normalized so that  $\alpha (u) \le \pi /2$. Then the condition 
 $d(z, A) > 2 $ guarantees that 
  $\alpha (\lambda z) \le \pi /3$.  However $ |z - \theta\lambda z| \ge
|z - \lambda z|$, then   
$$
\beta =\alpha (gz) \le  \alpha (\lambda z) \le \pi /3$$

\begin{figure}[htbp] %  figure placement: here, top, bottom, or page
   \centering
   \includegraphics[width=4in]{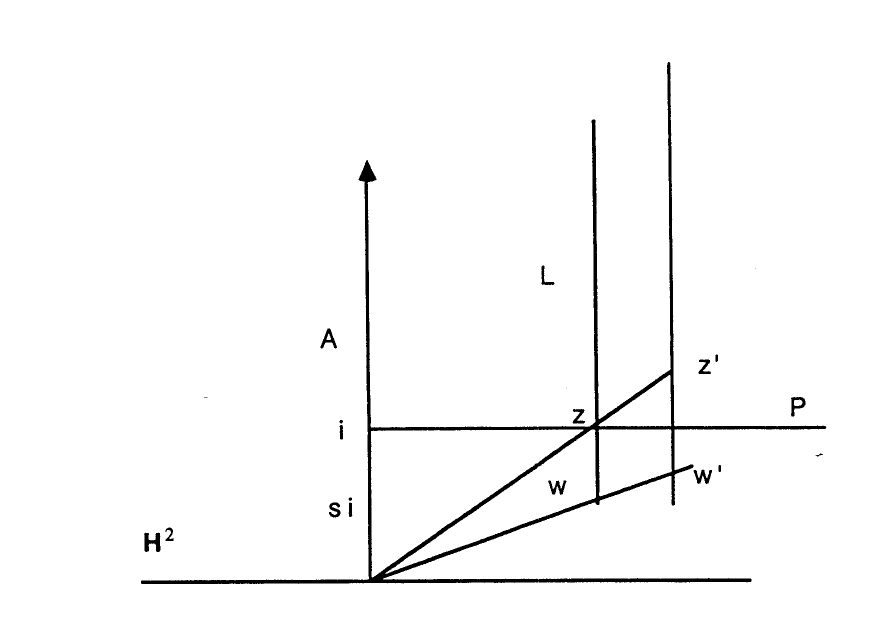} 
   \caption{}
   \label{F3}
\end{figure}

\medskip
Step 2. Due to the Step 1 it suffices  to consider the case:
$z, \ z' =\lambda z, \ w, \ w' = \lambda w \in \H ^2 \subset \C$ ,
 $ \arg (z) = \beta \le \pi /3 ; \ 
d(z , w) = d(z', w'), \ Re (z) = Re (w)$. (See Figure \ref{F3}.) Without loss of generality we can
 assume that $Im (z) = 1, y = Im (z')$.

Now we have:  $\rho = |z - z'|, y = \rho \sin ( \beta ) + 1 $,

$q = 2 \sinh  \frac{d(z, z')}{2} = \frac{\rho}{\sqrt{ y}} = 
 \frac{\rho} { \sqrt {1 + \rho \sin (\beta)}} $.

So $q^2 + \rho q^2 \sin (\beta) - \rho ^2 = 0$,
$$
\rho = (q^2 \sin (\beta ) + \sqrt {q^4 \sin ^2 (\beta ) + 4q^2})/2 
\ge q \ge 2 \sinh \frac{\nu}{2} \eqno(6)$$ 

On other hand we have: $\sin \beta \le \sqrt {3} /2$, so 
$\rho \le 2q^2 + 2 \le 8 \sinh ^2 (R/2) + 2 $,
$$ 
\sinh \nu /2 \le \sinh d(z, z'')/2 = |z - z''|/2 = \rho ( \cos \beta )/2 $$
$$
\le ( 8 \sinh ^2 (R/2) + 2) \cos \beta \le  8 \sinh ^2 (R/2) + 2 $$

 and $d(w, w') \ge d(w, w'') - d(w', w'') \ge d(w', w'') - R$.

Let  $s = Im(w)$, then $ d(z, w) = \log (1/s)$ and
$$
\sinh ( d(w, w'')/2) = |z - z''|/(2  s) \ge \rho /(4 s) 
\ge \sinh (\nu /2) /(2 s) \eqno(7) $$
$$ s \ge \frac{1}{2} \sinh ( \nu /2) \sinh  ^{-1} \frac{d(w, w'')}{2} \ge
\frac{1}{ 2} \sinh ( \nu /2) \sinh  ^{-1} (R) \eqno(8)$$
since $ d(w, w'') \le R +  d(w, w') $ and $ d(w, w') \le R $.
Now 
$$ d(z, w) = \log (1/s) \le \log (2 \sinh ^{-1} (\nu /2) \sinh R ) \eqno(9)$$

However $ \log \sinh a = (a^2 - 1)/2a $ and $2 \sinh b \le e^b $.
Therefore: $ d(z, w) \le R + \frac{ 1}{ \nu } $. 

 Lemma follows. \qed

\bigskip

Now suppose that 
 $a, b, z \in \H ^4 $ be points such that 
$d(z, [a, b]) \le R $. 

 Denote by $L_a , L_b $ the geodesic rays connecting the points $a, b $ and
$\infty $.

\bigskip

\begin{figure}[htbp] %  figure placement: here, top, bottom, or page
   \centering
   \includegraphics[width=4in]{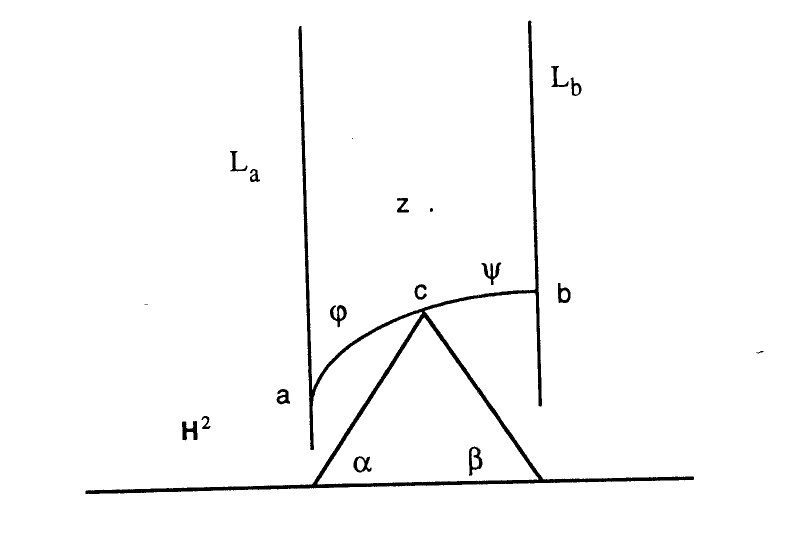} 
   \caption{}
   %\caption{example caption}
   \label{F4}
\end{figure}

\begin{prop}\label{p4}
%PROPOSITION 4. 
Under the conditions above we have:
$$
\min \{ d(z, L_a), d(z, L_b) \} \le 2 + R \eqno(10)$$
\end{prop}
\proof 
%PROOF. 
Denote by $c$ the point of $[a, b]$ such that $d(c, z) \le R $.
Let $L'_x $ be the geodesic containing $L_x$. Then we have:
$$
\cosh d(c, L'_a) = \sin^{-1}\alpha, \ \cosh d(c, L'_b) = \sin^{-1}\beta
\eqno(11)$$
(see Figure \ref{F4}) and $\alpha + \beta \ge \pi /2 $ so
$$
\sin^2\alpha + \sin^2\beta \ge 1 \eqno(12) $$

Now there are two possibilities (up to swapping $\beta \to \alpha $):

(i) $b_4 \ge w_4 $ for every $w \in [a, b]$

(ii) else.

Consider (ii). Then $\phi < \pi /2 \ , \psi < \pi /2 $ where
$\phi, \psi $ are nonzero angles of the triangle formed by
$L_a ,\ L_b ,\ [a, b]$. Therefore:
$$
d(c, L'_a ) = d(c, L_a) \ ,\ d(c, L'_b ) = d(c, L_b) \eqno(13)$$

Now if $\sin ^2 \alpha \le 1/4 $ then $\sin^{-1} \beta  \le 2 $.

This means that
$$
\min \{\cosh d(c, L_a) , \cosh d(c, L_b) \} \le 2 \eqno(14)$$ 

so $\min \{ d(z, L_a) , d(z, L_b) \} \le 2 + R $ in the case (ii).

Case (i). Then  $\phi < \pi /2 \ , \psi > \pi /2 $ however
$\alpha \ge \pi /4 $ (since the arc of the geodesic passing through $a, b$
is greater then the quarter of circle). 

Then $ 1/ \sin \alpha \le \sqrt{2} $ and $e^x / 2 \le \cosh x = \cosh (d(c,
L'_a )) \le \sqrt{2} $; $x \le \log (3) < 2 $.

However $\phi < \pi /2 $ , then $d(c , L'_a) = d(c, L_a) $ that means
$d(c , L_a) < 2 $.

Therefore 
$$
\min \{ d(z, L_a) , d(z, L_b) \} \le 2 + d(c, z) \le 2 + R \eqno(15)$$

\qed

\bigskip

2.5. 

%\ LEMMA 5. 

\begin{lem}\label{l5}
Suppose that $g$ is a loxodromic element with the axis
 $A$ as in Lemma \ref{l4};
 $a, b, z \in \H ^4 $ be points such that:
$d(z, [a, b]) \le R $,   $d(z, A) \ge 2 + R$,
 $\nu < d(g(z) , z) < R \ , \ d(a, g(a)) < R \ , \
d(b, g(b)) < R $, so the points $z, a, b$ have bounded index with respect
to $\langle g\rangle $.

Then $ \min \{d(z, a) , d(z, b) \} \le 4R + 6 + 1/k $,
where 
$$
k(R, \nu) = k = 2 (2+ R) \nu ^3 / \exp (18(2+ R) + 2 \nu ). $$
\end{lem}

\proof According to the Proposition \ref{p4} we can assume that 
$d(z, L_a) \le 2+ R $. Let $u \in L_a$ be a point such that
$d(z, u) < 2 + R$. Then we have: $d(u, A) > 2 $,
$ k < d(g(u), u) < 3R +4 $ (the last is by Lemma 3). Now we can apply
Lemma 4 to the $a, u$ to obtain: $d(a, u) < 3R + 4 = 1/k$ and
$d(z, a) < 4R + 6 + 1/k $.
 \qed

\begin{cor}
%COROLLARY 5 (
[The Key Property of Index] \label{c5}
Let $h$ be a parabolic or loxodromic
isometry of $\H ^4$,  $x \in \H ^4 $ is such that: $Ir_{\langle h \rangle} (x) > \mu $ and 
there is a geodesic segment  $L= [a, b]$ such that $d(x, L)
 < R $ for some  $R$  and
$$
\max \{\ind(x), \ \ind(a),\  \ind(b) \}  < R$$ 
 Then   either  
$$
\min \{ d(x,a), d(x,b)\} < C_+(R , \mu )= 4R +6 + 1/k$$
where $k=  2 (2+ R) \mu ^3 / \exp (18(2+ R) + 2 \mu ) $ (parabolic alternative)

 or $l(h) > C_-(R, \mu )= C_1(R + 2, \mu ) > 0 $
where the function $C_1$ is defined by Lemma 3
 (hyperbolic alternative).
\end{cor}
\proof 
%PROOF. 
Combine Lemma \ref{l3} and Lemma \ref{l5}.  \qed

Denote by $C(t, A) = \{w \in \H ^2 : d(w, A) = t\}$ the ``hypercycle" whose
 axis is the geodesic $A$.

\begin{prop}\label{p6}
%PROPOSITION 6. 
Let $z_1 , z $ belong to the same connected component of $C(t, A)$.
Then
 $$
d_{C(t, A)} (z_1 , z) \le 2\sinh (d(z_1 , z )/2) \eqno(16)$$
 where
$d_{C} $ is the metric on $C = C(t, A)$ induced from the hyperbolic plane.
\end{prop}
\proof We can assume that $|z_1| = 1, |z|= r$, $\log r $ is just the
distance between projections of $z_1 , z$ on the geodesic $A$.
 Then
$$
2\sinh (d(z_1 , z)/2) =
\frac{r-1}{ \sin (\theta) \sqrt {r} } \eqno(17)$$
 for $\cosh (t) \sin \theta = 1$. Here $\pi  - 2 \theta $ is the Euclidean
angle at the vertex of $C(t, A)$. 
Moreover, $a = d_C (z_1 , z) = \log (r) /\sin \theta$, $ a \sin \theta = \log
r $. Our aim is to show that:
$$
\log(r) \le \frac{r-1}{  \sqrt {r} } = \sqrt {r} - \frac{1}{\sqrt {r}}
\eqno(18)$$
if $r \ge 1$. Let $x = \sqrt {r}$, then $2\log x \le x - 1/x$ since for
$x = 1$ we have the equality and derivative of the left side is $\le$ of derivative of the right
side.
\qed 

%REMARK 5.
 \begin{rem}
 In the situation above we also have:
 $$
 \cosh t \le \frac{2}{\log r }\ \sinh \frac{d(z_1 , z)}{ 2 }. $$
\end{rem}

\bigskip

2.6. From now and until section 2.13, we let $q \in Isom(\H ^4)$ be a nonelliptic isometry.
 If $q$ is loxodromic then we
 shall suppose that $q$ fixes $\{ 0, \infty \}$ and $0$ is the repulsive
 point. If $q$ is parabolic then we shall assume that $q$ fixes the
point $ \infty $ (i.e. $q$ is an isometry of the Euclidean space).

Recall that $p: \H ^4 \to \H ^4 /\langle q\rangle$ is the covering map;
$$
{\mathcal K} (\langle q\rangle, \nu) = \{ x \in \H ^4 : \inf \{ d(x, g(x)): g \in \langle q\rangle \setminus 
\{ 1\} \} \le \nu \} \ ,$$
${\mathcal T} (\langle q\rangle, \nu) = p({\mathcal K} (\langle q\rangle, \nu))$.

For the element $T \in \langle q\rangle $ we denote by $T_{\theta }$ its rotational
component and $T_{\lambda } = T_{\theta }^{-1}T$. Also let 
 $0 \le \theta _T \le \pi$ be the
angle of rotation of $T_{\theta }$; denote by $\lambda _T $  either the
 Euclidean
distance $|X - T_{\lambda }(X)|$ (in the case of parabolic $T$) or the 
coefficient of similarity (in the loxodromic case). If 
the rotational component of $q$ is non-trivial, let $R_x$  be the Euclidean
distance from the point $x \in \H ^4$ to the Euclidean plane of rotation
$L_q$. For loxodromic $q$ define $A_q$ to be the geodesic axis of $q$.
If $\theta _q = 0$ then we set $L_q = \emptyset , R_x = 0$.
Direct calculations show that: 
$$
|T(x) - x|^2 = 2R_x^2 (1- \cos \theta _T) + \lambda _T^2 \eqno(19)$$ 
(in the parabolic case)
$$
|T(x) - x|^2 = (2R_x^2 (1- \cos \theta _T) + \lambda _T^2)|x|^2 \eqno(20)$$
(in the loxodromic case).
Therefore 
$$
  2 \sinh ^2 (d(x, T(x)/2)=  (2R_x^2 (1- \cos \theta _T) +
 \lambda _T^2)/x_4^2 \eqno(21)$$
 (in the parabolic case)
$$
 2 \sinh ^2 (d(x, T(x)/2)=  (2R_x^2 (1- \cos \theta _T) +
 (\lambda _T - 1)^2/\lambda _T )|x|^2/ x_4^2 \eqno(22)$$
(in the loxodromic case).

We shall denote $\sqrt{2} \sinh (d(T(x),x)/2)$ by $x_T$; set    
$$
 x^\theta _T = |x|R_x \sqrt {2(1- \cos \theta _T)} /x_4 $$
 and 
$x^\lambda _T = \sqrt {x_T^2 - (x^\theta _T)^2 }$.

The formulas (21), (22)  imply that the domain $ {\mathcal K} (\langle g\rangle , \nu)$ is
 convex
near its smooth boundary points, where the boundary from the Euclidean
point of view is a piece of a cone (in loxodromic case) or cylinder
 (in parabolic case) over an ellipsoid of
revolution or 2-sheeted hyperboloid of revolution. In smooth points the
boundary  $\partial {\mathcal K} (\langle q\rangle, \nu)$ is given by the equation
$$
x_{q^{n(x)}} = \sqrt {2} \sinh (\nu /2) \eqno(23)$$
 However in nonsmooth boundary
points the domain  $ {\mathcal K} (\langle q\rangle, \nu)$ is not locally convex. Nevertheless
the following remark will be important to us.

\begin{figure}[htbp] %  figure placement: here, top, bottom, or page
   \centering
   \includegraphics[width=4in]{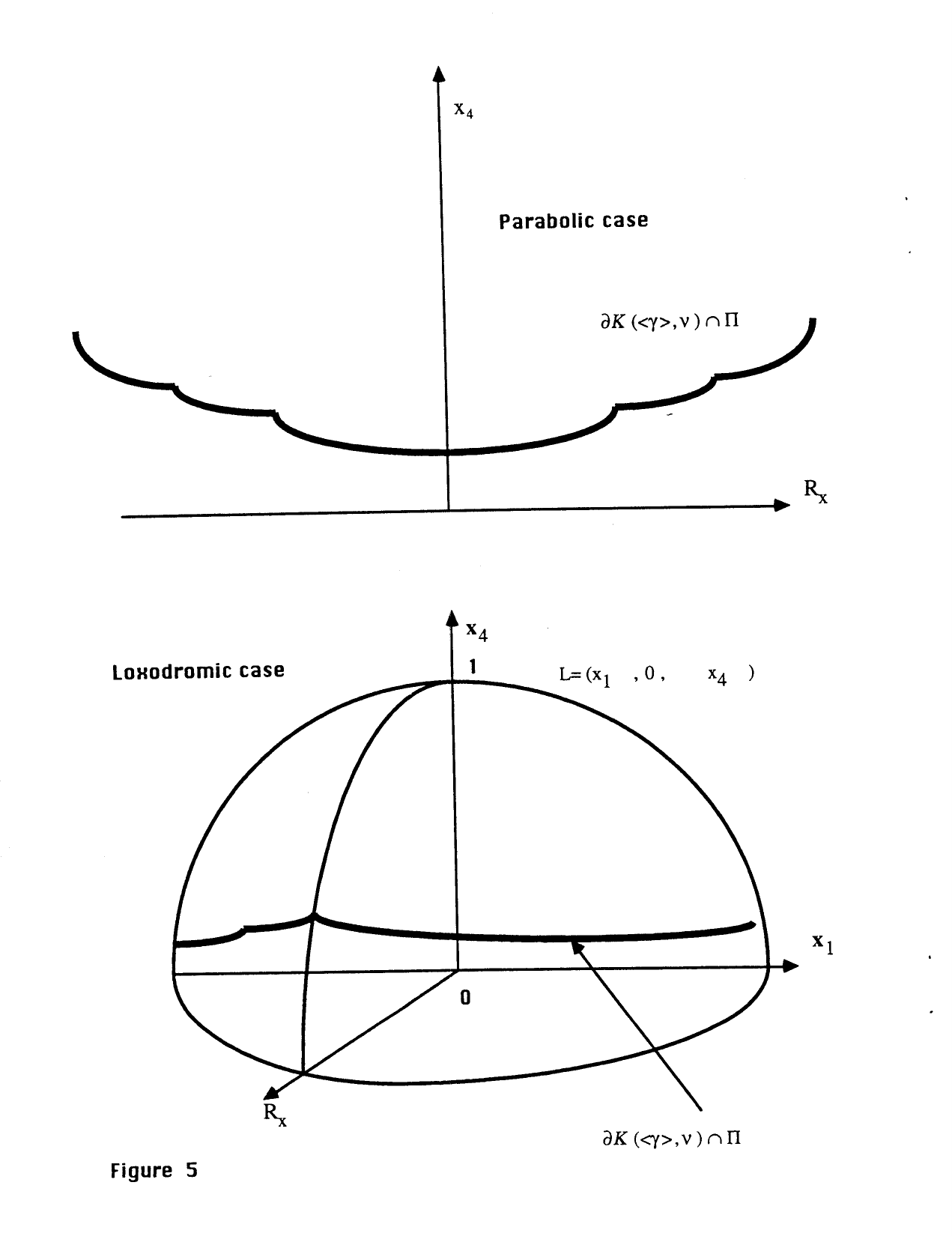} 
  \caption{}
   \label{F5}
\end{figure}

%REMARK 6. 
\begin{rem}\label{r6}
Take an arbitrary fiber  $\Pi _t$ of the canonical
fibration and let  $\H ^2 \subset \Pi _t$ be any hyperbolic plane which
 contains the axis of rotation $L_q \cap \Pi _t$. Then the curve 
$$
\partial {\mathcal K} (\langle q\rangle, \nu) \cap  \H^2 \cap \Pi _t \eqno(24)$$ 
is given by the equation $x_4 = x_4 (R_x)$ which is an increasing function
on $R_x$. In the both parabolic and loxodromic cases this equation is the
equation of hyperbola (in smooth points); the difference between parabolic
and loxodromic cases is that in the last case the domain of the function
is bounded by (say) $R_x \le 1$. (See Figure \ref{F5}.)
\end{rem}

Suppose $H \subset Isom(\H ^4)$ is a  virtually  abelian discrete subgroup;
$H(\infty )
= \infty$. 
If $H$ contains a loxodromic element $q$ then denote by $ \H ^4_*$ the
 complement in  $ \H ^4$ to the  axis $A= A_q$ of $q$
  and we set $\H ^4_* = \H ^4 $  in the parabolic case.

Let $\varphi = \varphi _{H, \nu} : a \in  \H ^4 _* \to \partial {\mathcal K} 
(H, \nu)$
 be the projection
$
\varphi (a) = $ the point of intersection of $\partial {\mathcal K}(H, \nu) $
with the geodesic
containing $ a \ , \infty $ \ (in the parabolic case),
$
\varphi (a) = $ the point of intersection of $ \partial {\mathcal K}(H, \nu) $
with the geodesic ray 
containing $ a $ and orthogonal to  $A_q$ (in the loxodromic case).

\medskip 
2.7. 

\begin{lem}
%LEMMA 6. 
[Compare \cite{B}, \cite{SY}, \cite{HI}] \label{l6}
The map  $\varphi $ is well-defined.
\end{lem}
\proof  Firstly, formulae (21-22) imply  that the set ${\mathcal K}(H, \nu) $ is
 star-like with
respect to the point $\infty $. 
If $H$ contains  loxodromic element $q$  and $z \in \partial {\mathcal K}(H, \nu)$
 then we can use the fact that the whole Euclidean ray
$K_z = \{ c \cdot z : c \in \R _+ \}$
is also contained in $\partial {\mathcal K}(H, \nu)$. Therefore the domain 
in the Euclidean plane between the rays $A_q$ and $K_z$ is contained
in $ {\mathcal K}(H, \nu)$. Hence the point of intersection defining 
 $\varphi $ in the loxodromic case is unique. 
 Now it follows from (3-4) that the intersection of 
 ${\mathcal K}(H, \nu) $  with any geodesic asymptotic to the point 
$\infty $ is non-empty.
 \qed

\begin{prop}\label{p7}
%PROPOSITION 7. 
Suppose that $\ind_g(a) \le R$ ; $\nu \le R$. (i) Then either
$$
\min \{ d(a, \varphi (a))| a \in {\mathcal K}(\langle g\rangle , \nu ))\} 
\le C_+(R, \nu ) $$ 
or
$l(g) \ge C_{-}(R, \nu ) $ and 
$$
\cosh d(a, A_g) \le  \frac{2\sinh R/2}{  C_{-}(R, \nu )}.  
$$
 (ii) If $g$ is parabolic, then 
$$
\min \{ d(a, \varphi (a))| a \in 
{\mathcal K}(\langle g\rangle , \nu )) \} \le  1 + R/2 - \nu /2 .$$
\end{prop}
\proof (i) Follows from Lemma \ref{l5} and Corollary \ref{c5}. (ii) Suppose that 
$d(b, g(b)) \ge \nu $. Then 
$$
\frac{\sinh ^2 (R /2)}{ \sinh ^2 (\nu /2)}
\ge b_4^2/a_4^2,
$$
 where $b= \varphi (a)$. Now the statement (ii)
follows from  direct calculations. \qed

\bigskip
2.8. Let $x, z \in \H ^4$ be a pair of distinct points, $T \in\langle q\rangle  \setminus \{ 1 \}$.
Define the film $S = S_{Txz}$ in $\H ^4$ connecting points $x , z , T(z), 
 T(x) $ as follows.

Let $T = T_{ \theta} \circ T_{ \lambda} ;  T_{ \theta} = \exp (\xi ),
 T_{ \lambda} = \exp (\zeta ) ; \ \xi = \xi _T , \zeta = \zeta _T \in so(4,1)$,
where 
$$\exp : {t \in [0, 1]}\cdot \xi \mapsto SO(4,1)$$ 
is injective.

 Then we set 
$$
S^\lambda = \bigcup _{t \in [0, 1]}  \exp(t \zeta  )([x, z]) \eqno(25)$$
$$
S^\theta = \bigcup _{t \in [0, 1]}  \exp(t \xi )(g_{\lambda }[x, z]) \eqno(26)$$
$$
S = S^\lambda \cup S^\theta $$
$$
 \partial_x S = \bigcup _{t \in [0, 1]}  \exp(t \zeta  )(x) \cup
\bigcup _{t \in [0, 1]}  \exp(t \xi )(T_{\lambda } x) \eqno(27)$$
$$
 \partial_z S = \bigcup _{t \in [0, 1]}  \exp(t \zeta  )(z) \cup
\bigcup _{t \in [0, 1]}  \exp(t \xi )(T_{\lambda } z) \eqno(28) $$
$$
\delta S =  \partial_x S \cup  \partial_z S \eqno(29)$$

 The films $S = S_{Txz}$ constructed above will be called ``ruled films". 
It's easy to  see that number of points of 
transversal intersection of ruled film with another ruled film (or geodesic
plane in $\H ^4$ ) is at most $8$.

%REMARK 7. 
\begin{rem}\label{r7}
The film $ S_{Txz}$ is contained in $ {\mathcal K}(\langle g\rangle , \nu)$.
\end{rem}

The ruled film $ S_{Txz}$ is said to be {\em in general position} if:

(i)  $ [x, z] \subset \H_* ^4 - L_q$;

(ii) In the loxodromic case let $\H ^3_{[x, z]}$ be the geodesic hyperplane
 of $\H ^4$ which contains
$[x, z] \cup A_q$. Then $\H ^3_{[x, z]}$ is not  orthogonal to $L_q$.

\medskip 
2.9. Suppose that we are given a ruled film $S_{gxz}$,
then $ \partial_x S_{gxz}$ consists of two arcs:
$$
\delta _x ^\lambda = \bigcup _{t \in [0, 1]}  \exp(t \zeta  )(x) $$
$$
\delta _x ^\theta =
\bigcup _{t \in [0, 1]}  \exp(t \xi )(T_{\lambda } x) .$$

Then take totally geodesic regions $D(x , g_\lambda (x))$ and 
$D(g_\lambda (x) , g (x))$ bounded
by $\delta _x ^\lambda \cup [x, g_{\lambda } x]$ and $\delta _x^\theta
 \cup [ g_{\lambda } x \ , g x]$ respectively.

 Define the totally geodesic
regions $D(z , g_\lambda (z))$ and $D(g_\lambda (z) , g (z))$ analogously.

\bigskip

%DEFINITION 1. 

\begin{figure}[htbp] %  figure placement: here, top, bottom, or page
   \centering
   \includegraphics[width=4in]{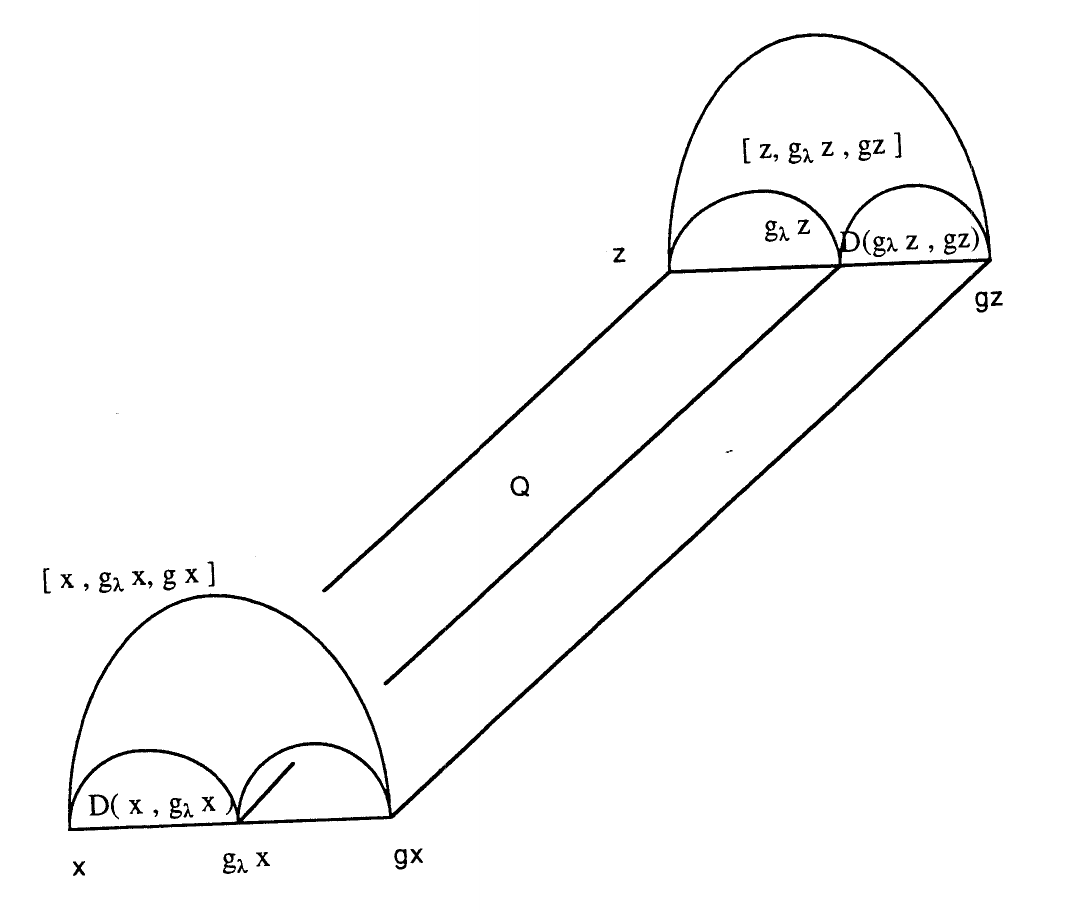} 
   \caption{}
  % \caption{example caption}
   \label{F6}
\end{figure}

\begin{defn}\label{d1}
The film 
$$
{\tilde Q}_{gxz} = S_{gxz} \cup D(x , g_\lambda (x)) \cup
 D(g_\lambda (x) , g (x))
 \cup $$
$$
D(z , g_\lambda (z)) \cup D(g_\lambda (z) , g (z)) 
\cup [x,  g_{\lambda } x \ , g x] \cup  [z,  g_{\lambda } z \ , g z]
$$
is called the {\em extended ruled film}. If $g \in G \subset Isom(\H ^4)$
then the projection $Q_{gxz}$ of ${\tilde Q}_{gxz}$ to $\H ^4 /G$ is called 
 the {\em extended ruled annulus}.
Define
$$
  \partial_{ x} Q_{gxz} ,   \partial_{ z} Q_{gxz}$$
to be the projections of $[x, gx], [z, gz]$ in  $\H ^4 /G$. Then 
$$
  \partial_{ x} Q_{gxz} \cup   \partial_{ z} Q_{gxz} = \partial Q_{gxz}$$
(see Figure \ref{F6}).
\end{defn}

\medskip
2.10. 

\begin{lem}\label{l7}
%LEMMA 7. 
Suppose that the  film $S_{gxz}$ is in general position. Then
 $p(  \partial_x S_{gxz} )$, $p(  \partial_z S_{gxz} )$ are homologous
in 
$$
\partial {\mathcal T}(\langle g\rangle , \nu) - p(L_q \cap \partial {\mathcal K}(\langle g\rangle , \nu))\ .$$
Moreover, the intersection
 $p( S_{gxz}) \cap  p(L_q \cap \partial {\mathcal K}(\langle g\rangle , \nu))$ is empty.
\end{lem}
%(a) $p( \partial  S_{gxz} ) = p( \partial  S_{gxz'} )$;
%(b) $p \varphi _q ( S_{gxz'})  \cap  p(L_q \cap \partial {\mathcal K}(\langle g\rangle , \nu))
%= \emptyset$.

\proof By the property (i) of film in general position we have:
  $ [x, z] \cup  L_q = \emptyset $. However $S_{gxz}$ results  
via moving the segment  $[x, z]$ by elements contained in $Z_{\langle g\rangle }$
 which leave $L_q$ invariant.
%If $\theta _q =0$ then there is nothing to prove. So
%we can use the fact that $\theta _q \ne 0$ to find a point $z' =
%q^k(z)$ such that $[z' , x] \cap L_q = \emptyset$. Now the film  $S_{gxz'}$
%is obtained by moving the segment $[z' , x]$ by elements of $Z_{\langle g\rangle }$.
%Therefore, all conditions are fulfilled. 
\ \ \ \ \qed
 
 \medskip 
2.11. Suppose that $x, y, z, w \in  \partial {\mathcal K}(\langle g\rangle ,
 \nu)) \cap \Pi _t$ --- the same fiber of the canonical fibration
associated with $q$; $g, h \in\langle q\rangle \setminus \{ 1\}$,  
$$
\max \{ d(g(x), x), d(g(y), y), d(g(z), z), d(g(w), w) \} \le C\ , \eqno(30)$$
 $S_{gxz} , S_{hyw}$ are in general position,  
$p(\partial  S_{gxz}) \cap  p(\partial  S_{hyw}) = \emptyset $
and the intersection $p(\partial  S_{gxz}) \cap  p(\partial  S_{hyw})$ is 
transversal.
Our purpose is to estimate the algebraic intersection number between the 
annuli
$$
p(S_{gxz}, \partial S_{gxz} ) , p(S_{hyw},  \partial S_{hyw}) 
\subset ({\mathcal T}(\langle q\rangle, \nu), \partial {\mathcal T}(\langle q\rangle,
 \nu))$$
or, equivalently, the algebraic intersection number between  
$S_{gxz} $ and $\langle q\rangle (S_{hyw})$. Notice that this number depends only
on
$$
p( \partial S_{gxz} ) , p( \partial S_{hyw}) 
\subset  \partial {\mathcal T}(\langle q\rangle, \nu)$$
and does not depend on the relative cycles in $ {\mathcal T}(\langle q\rangle, \nu)$ with
the boundaries $p( \partial S_{gxz} )$, $p( \partial S_{hyw})$.

\begin{thm}\label{t4}
%THEOREM 4. 
Under conditions above, we have the following estimate for the 
algebraic intersection number:
$$
|<p(S_{gxz}),  p(S_{hyw}) >| \le N(C, \nu ) \eqno(31)$$
where 
$$
N(C, \nu ) = ( \exp ^3 (4C +4))/\nu ^3 \eqno(32)$$
\end{thm}
\proof  First, consider the most interesting case
$\theta _q \ne 0$. Lemma 7 implies that if $p(Z_{\langle q\rangle}(x))$ does not divide 
say $p(Z_{\langle q\rangle}(y))$ 
from $p(Z_{\langle q\rangle}(w))$
then the annulus ${\mathcal A}= p(S_{hyw}) $ can be deformed
 $rel ( \partial {\mathcal A})$ to
a new annulus ${\mathcal A} '$ so that ${\mathcal A}' \cap
 p(  \partial_x S_{gxz}) = \emptyset $ and 
 $$
  \#({\mathcal A}' \cap
 p(  \partial_z S_{gxr})) \le \#({\mathcal A} \cap
 p(  \partial_z S_{gxz})).$$

Therefore (up to a change of notation) our problem is reduced to the case 
$$
R_x \le R_y \le R_z \le R_w \eqno(33)$$
(if $q$ is parabolic) and
$$
\min \{R_x ,  R_z \} \le R_y \ \ , \ \ \min \{R_y ,  R_w \} \le R_z
\eqno(34)$$
(if $q$ is loxodromic).
The monotonicity of the boundary $\partial {\mathcal K}(\langle g\rangle, \nu)$ (see Remark 6)
implies that:
$$
\min \{x_4 ,  z_4 \} \le y_4 \ \ , \ \ \min \{y_4 ,  w_4 \} \le z_4
\eqno(35)$$
If either 
$$
 \diam ( \varphi (S_{hyw}) \cap (\langle q\rangle  \partial_z S_{gxz} ))\le Const(C,  \nu)
= 2\sqrt {C+1}
  \eqno(36a)$$
or 
$$
\diam (\varphi (S_{gxz}) \cap (\langle q\rangle  \partial_y S_{hyw} )) \le Const(C,  \nu)
= 2\sqrt {C+1}  \eqno(36b)$$
then we can apply Lemma 3  to obtain
$$
|<p(S_{gxz}),  p(S_{hyw} )>| \le \exp ^3 (2Const(C, \nu ))/\nu ^3 =
\exp ^3(4C +4)/\nu ^3$$
 So, our goal is to obtain one of such estimates. The inequalities (35) imply that
$$
y^\lambda _g \le \max \{ x^\lambda_g , z^\lambda_g \} \le \sqrt {2} \sinh (C/2) ;
z^\lambda _g \le \max \{ y^\lambda _g , w^\lambda _g \} \le \sqrt {2} \sinh (C/2)$$
 Therefore, consider 
$$
z^\theta _g/z^\theta _h = y^\theta _g/y^\theta _h = \frac{\sin \theta _g/2}{\sin \theta
_h/2}= \alpha. $$

 Now if $\alpha \ge 1$ then  $z^\theta _g \ge z^\theta _h$; if  $\alpha \le 1$ then 
 $y^\theta _g \le
y^\theta _h$ . So either
$$
\ind_h(z)= d(z, h(z)) \le \sqrt {\mathrm{arcsinh} (4\sinh ^2 C/2)} = C' \le  
\sqrt {C+1} \eqno(37)$$ or
$$
\ind_g(y)= d(y, g(y)) \le \sqrt {\mathrm{arcsinh} (4\sinh ^2 C/2)} = C' \le \sqrt {C+1}\eqno(38)$$
We shall assume that (37) holds.

Notice that the intersection
$$
\varphi (S_{hyw}) \cap (\langle q\rangle  \partial_z S_{gxz} )$$
 is contained in
$\varphi (S_{huw}) \cap Z_{\langle g\rangle }( z )$. Moreover, because our films are in
general position (condition (ii)) we have
$\varphi [u, w] \cap  Z_{\langle g\rangle }( z ) = \{s_1 , s_2 \}$ where $s_1 , s_2 $ can
coincide. 

Therefore
$$
\varphi (S_{huw}) \cap Z_{\langle g\rangle }( z ) =$$ 
$$
\bigcup _{i= 1, 2} \bigcup _{t \in [0, 1]} \exp(t \zeta _h  )(s_i) \cup
\bigcup _{t \in [0, 1]} \exp(t \xi _h)(h_{\lambda  } s_i . \eqno(39)$$
However $\ind _h(s) = \ind_ h(z) \le C'$ due to (37); hence the diameter
 of the set in (39) 
is bounded from above by $2C'$  and we are done. 

\medskip 
2.12. Now suppose that the rotational component of $q$ is trivial. Then the
condition (30)  implies that $m, n \le [C/\nu]$. Take 
 $x, y, z, w \in \Pi _t$.  So the films $ S_{hyw} , S_{gxz}$
lie in the union of  $[C/\nu]$ images of a convex fundamental domain
of the group $\langle q\rangle$. Hence the intersection number is at most 
$8 [C/\nu] \le N(C, \nu )$.

\qed

\medskip 
2.13. Consider the case when $g, h$ are parabolic elements which belong to
a  virtually abelian group $\Gamma \subset Isom(\H ^4)$ such that $\Gamma$ is not 
a cyclic group. Again suppose that $x, y, z, w \in 
  \partial {\mathcal K}(\Gamma ,
 \nu)$  are points such that the condition (30) is fulfilled. 
Denote by $p: \H ^4 \to \H ^4 /\Gamma \ $ the universal covering. Our aim is
estimate the intersection number between $pS_{gxz}$ and $pS_{hyw}$ in terms
of $C, \nu$.

\begin{thm}\label{t5}
%THEOREM 5.  
The intersection number between $pS_{gxz}$ and $pS_{hyw}$
is bounded from above by the constant 
$$
N'(C, \nu )= (\exp (  9\nu + 6))/\nu ^3 + 96C/\nu + 120000\exp (72C)/\nu ^3 
\eqno(40)$$
\end{thm}
%96C/\nu + 40\exp ^3(24C)/\nu ^3$$

\proof 2.14. Let $\Gamma _0 \subset \Gamma $ be a  maximal abelian
subgroup in $\Gamma$; $|\Gamma : \Gamma _0| \le 12$. Denote
by $g_0 = g^{n_g}$, $h_0 = h^{n_h}$ generators of the groups $\langle g\rangle  \cap \Gamma _0$, 
$\langle h \rangle \cap \Gamma _0$ respectively. Then ${n_g}$ and ${n_h}$
are not greater than $12$. Observe also that 
$$
 \partial {\mathcal K}(\Gamma _0,
 \nu) \subset  {\mathcal K}(\Gamma , \nu) -  {\mathcal K}(\Gamma , \nu /12) $$
and $\partial {\mathcal K}(\Gamma _0, \nu )$ is a horosphere.
For every $u \in \partial {\mathcal K}(\Gamma , \nu )$ by (21) we have
$$
d(u, \partial {\mathcal K}(\Gamma , \nu /12)) \le \log(\sinh(\nu /2)/\sinh(\nu
/24))
\le \nu /2 + 1 .\eqno(41)$$ 
Therefore for $u$ as above
$$
d(u, \partial {\mathcal K}(\Gamma _0 ,\nu )) \le \nu /2 .\eqno(42)$$

Denote by $ p_0:  \H ^4 \to \H ^4 /\Gamma _0$ 
the universal cover.

 2.15. First we reduce the problem to the abelian subgroup $\Gamma _0$.

(i) $|<pS_{gxz} , pS_{hyw}>| \le 12|<p_0(\langle g\rangle S_{gxz} ), p_0(\langle h \rangle S_{hyw})>|$,
where $HB$ denotes the orbit of the set $B$ under the group $H$;
$$
p_0(\langle g\rangle S_{gxz} )= p_0(S^0_{gxz}); (S^0_{gxz})= S_{gxz} \cup gS_{gxz} \cup g^2S_{gxz} ... \cup
g_0S_{gxz}$$

(ii) $\diam(\{u, Tu, T^2u, ..., T_0u \}) \le 12C $ for $u \in \{x, z \},
T= g$ and $u \in \{ y, w \}, T= h$.

\begin{figure}[htbp] %  figure placement: here, top, bottom, or page
   \centering
   \includegraphics[width=4in]{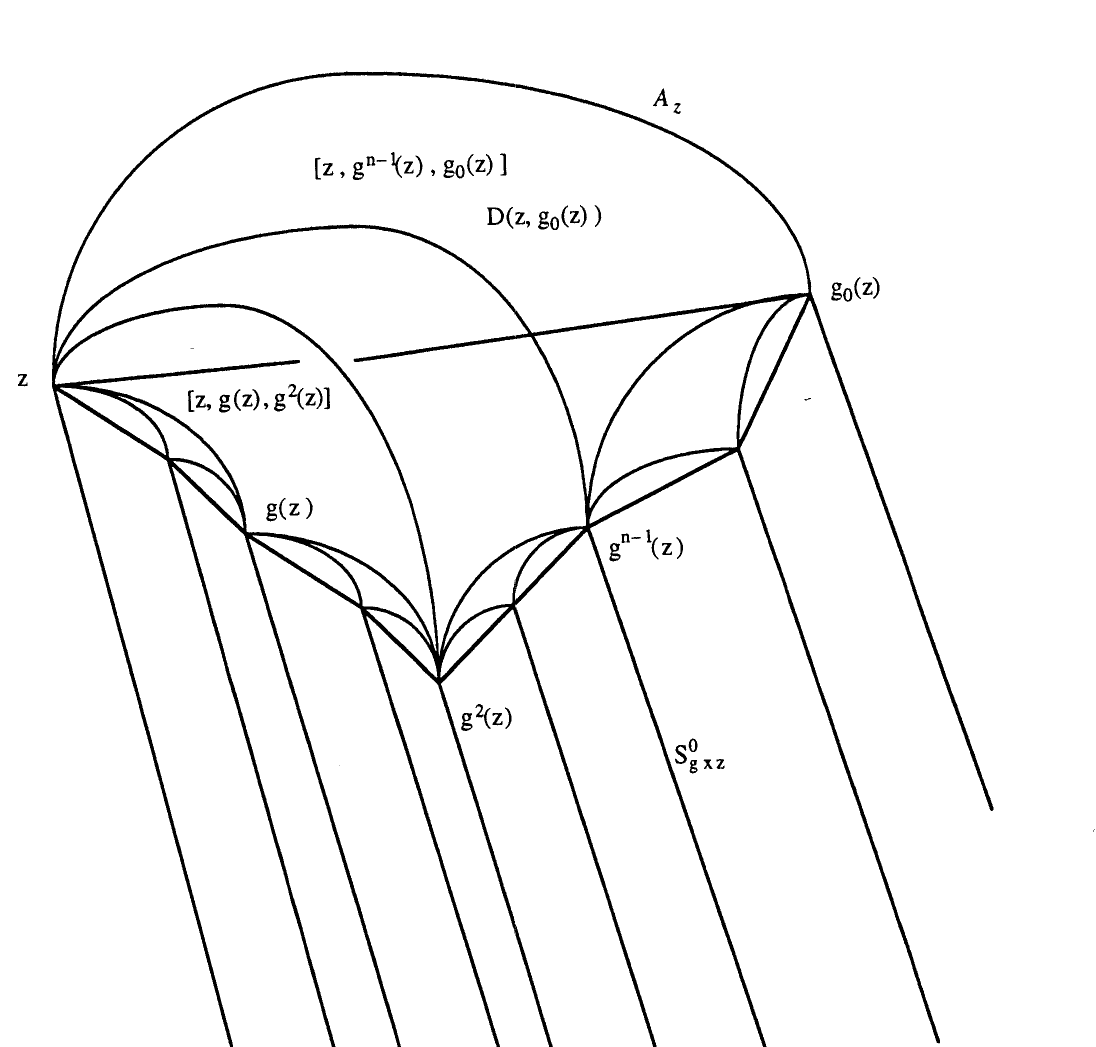} 
   \caption{}
  % \caption{example caption}
   \label{F7}
\end{figure}

Now we enlarge the films $S^0_{gxz}, S^0_{hyw}$ in the following way. Let 
$u \in \{ x, y, z, w \}$, \ $T \in \{ g, h \}$ be the corresponding
transformation. Then take the union
$$
{\mathcal A}^-_u = [u, Tu, T^2u] \cup [u, T^2u, T^3u] \cup ... \cup
 [u, T^{n_u-1}u, T^{n_u}u]$$
$$
{\mathcal A}_u = D(u, T_\lambda (u)) \cup D( T_\theta (u), T(u)) \cup 
D(T(u), T_\lambda (Tu)) \cup 
 ... $$
$$
\cup D(T^{n_u-1}u), g_\lambda ( T^{n_u-1}u)) \cup  
D( g_\lambda ( T^{n_u-1}u), T_0(u))
\cup D(u, T^{n_u}u) 
\cup {\mathcal A}^-_u$$
(see 2.9). Then set 
$$
S^+_{gxz}= S^0_{gxz} \cup {\mathcal A}_x \cup {\mathcal A}_z$$
$$
S^+_{hyw}= S^0_{hyw} \cup {\mathcal A}_y \cup {\mathcal A}_w$$ 
(Figure \ref{F7}).

All extra pieces that are attached to $ S^0_{hyw},  S^0_{gxz}$ have diameter
$\le 12C$ and their number  is bounded by $100$.

 Therefore 
$$|<p(S_{gxz} ), p(S_{hyw})>| \le 12|<p_0(S^+_{gxz} ), p_0(S^+_{hyw})>| + 
12\cdot 10000\exp^3(24C)/\nu ^3$$
$$
|<p_0(S^+_{gxz} ), p_0(S^+_{hyw})>| = 
|<p_0(S_{g_0xz}) ,
p_0(S_{h_0yw})>|$$
since $\partial  p_0(S^+_{gxz} ) = \partial p_0(S_{g_0xz} )$, $\partial 
 p_0(S^+_{hyw} ) = \partial p_0(S_{h_0yw} )$.

%+ (3I+2)(\exp^3(2 C \cdot I + 2\nu))/\nu ^3$ (by Lemma 3).

Denote by $\Gamma _1$ the maximal subgroup in $\Gamma _0$ which has the rank 2
and contains $g_0 , h_0$. We have two possibilities:

(a) $\Gamma _0$ has rank 3. 
(b) $\Gamma _0$ has rank 2.

In the case (a) denote by $k_0 \in \Gamma _ 0$ an element such that 
$\Gamma _ 0 = <k_0> \oplus \Gamma _1 $.

Choose a fundamental domain $\Phi$ for action 
of $<k _0>$ in  the horosphere
$\partial {\mathcal K}(\Gamma _0 \nu )$ such that $\Phi $ is bounded by a pair of 
Euclidean hyperplanes in $\partial {\mathcal K}(\Gamma _0 ,\nu )$. In the case
(b) let $k_0 = 1$ and $\Phi = \partial {\mathcal K}(\Gamma _0 ,\nu )$.

Let $\varphi : \H ^4 \to \partial {\mathcal K}(\Gamma _0 ,\nu )$ be the projection
defined in 2.6. Then $d(u, u'= \varphi u) \le 1 + \nu /2$ for every $u \in
\partial {\mathcal K}(\Gamma , \nu )$.

For $u \in \{ x, z, y, w \}$ we choose $T_u \in <k_0>$ such as
$T_uu' \in \Phi$. Denote $T_u(u')$ by $u''$. Then above we substitute the films
$S_{g_0xz} , S_{h_0yw}$ by  $S_{g_0T_x(x')T_z(z')} , S_{h_0T_y(y')T_w(w')}$.

   The intersection number between projections of the ruled films $S_j$ in
$\H ^4/\Gamma _0$ depends only on projections of their boundaries
$\delta S_j$. Therefore: 
$$
|<p_0(S_{g_0xz}) , p_0(S_{h_0yw})>| \le (\exp (6+  9\nu))/\nu ^3
+ |<p_0(S_{g_0x''z''}) , p_0(S_{h_0y''w''})>| .\eqno(43)$$

Next notice that 
 $\Gamma _1 (\varphi S_{g_0x''z''}), \Gamma _1(\varphi S_{h_0y''w''}) \subset 
\Phi $
are Euclidean parallelograms.
Hence we can estimate the intersection number 
$|<p_0(S_{g_0x''z''}) , p_0(S_{h_0y''w''})>|$ in the same way as in 2.11 by 
 looking at the intersections of the  $\Gamma _1$-orbit of $ \varphi \delta S_{g_0x''z''}$
 with $ \varphi S_{h_0y''w''}$.
However
 $$
|\Gamma _1 : <g_0 > \oplus < h_0>| \le 12C/\nu \eqno(44)$$
which implies that 
$$
|<p_0(S_{g_0x''z''}) , p_0(S_{h_0y''w''})>| \le 96C/\nu \eqno(45)$$

Thus, the final estimate is 
$$
<pS_{gxz} , pS_{hyw} > \le 
(\exp (6+  9\nu))/\nu ^3 + 96C/\nu + 120000\exp (72C)/\nu ^3 \eqno(46)$$
\qed

\section {Proof of  Theorem \ref{t1}  under the maximality condition}

\bigskip
3.1. {\bf Notations.} 

Fix  some Margulis constant $\mu$ for the 
4-dimensional hyperbolic space. So if $h_1 \ , h_2 \in Isom(\H ^4)$
generate a discrete group $H$ and
$d(x, h_i(x)) \le \mu \ (i= 1, 2)$ for some $x \in \H ^4$ then  $H$ is a 
  virtually
abelian group.
Recall that the hyperbolic 4-space is realized as the upper half-space
 $\R ^4 _+$.
 By $p:\H ^4 \to M$ we shall denote the
universal covering of $M$; its deck-transformation group is $G$.
Let $S$  be a triangulated Riemannian surface so that the edges of
triangulation are geodesic arcs.
 Then the map
$$
f: S \to M$$
is called piecewise-geodesic (p-g) if the restriction of $f$ on every simplex
is a totally geodesic map. 
By $M_{(0, \mu]}$ and $M_{ [\mu , \infty )}$ we denote $\mu $-thin and
$\mu $-thick parts of $M$ respectively. The components of  $M_{(0, \mu]}$ 
are 
Margulis tubes (see definitions in Section 2). 
If $\Delta $ is a triangle then we shall denote by $\dot \Delta $ the set of
vertices of $\Delta $. 
By $l_N(\gamma )$ we  denote the length of the curve $\gamma $ in the
metric space $N$. If a transformation  $h \in Isom(\H ^4)$ is parabolic 
then we set $A_h = \emptyset$; if $h$ is loxodromic then $A_h$ is the axis of
$h$ (invariant geodesic).

\begin{figure}[htbp] %  figure placement: here, top, bottom, or page
   \centering
   \includegraphics[width=4in]{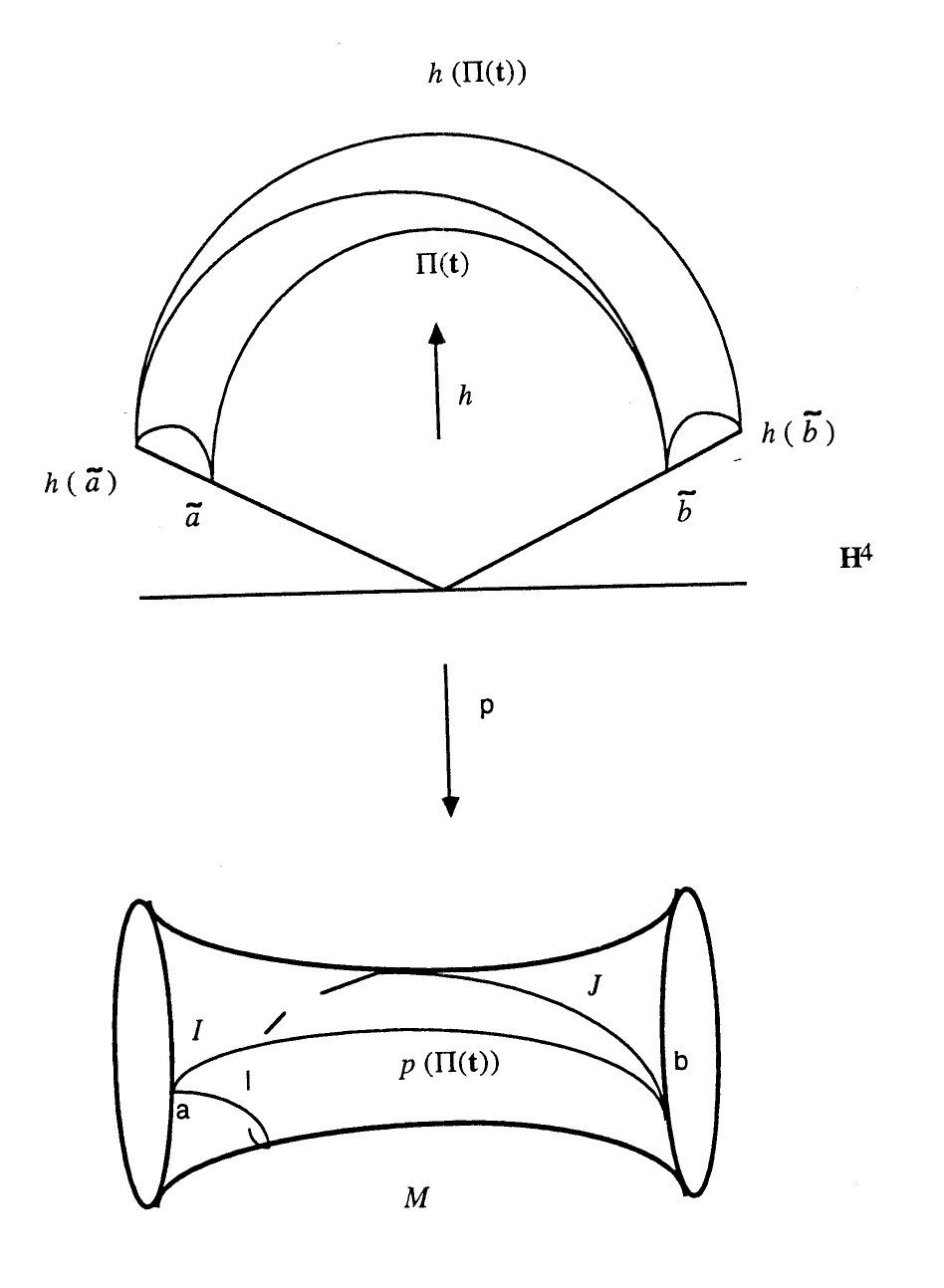} 
   \caption{}
   %\caption{example caption}
   \label{F8}
\end{figure}

\begin{defn}\label{d2}
%DEFINITION 2. 
Suppose $H \in Isom (\H ^4)$ is a  virtually abelian discrete
group, $h \in H - \{1\}$; $a, b \in {\mathcal T}(H, \nu )$ belong to one and the
same fiber $\Pi _t$ of the canonical foliation associated with $h$. Then
we define {\em the p-g annulus }$F_{hab}$ as follows.  First connect $a, b$
 by the geodesic segments $I,
J$
so that  $I \subset \Pi _t $, the closed loop $I \cup J$
is homotopic to $h$. Denote by $\gamma _a , \gamma _b$ the shortest loops
in $\H ^4/H$ which contain $x, y$ and homotopic to $h$. Then take the pair
 of geodesic triangles
in $M$ whose edges are $I, J , \gamma _a$ and $I, J , \gamma _b$ respectively.
The union of these triangles is the desired p-g annulus $F_{hab}$. 
(See Figure \ref{F8} for a lift of $F_{hab} $ in  $\H ^4$.) 
\end{defn}

 %REMARK 8. 
 \begin{rem}\label{r8}
 In general, $F_{hab}$ is not entirely contained in
 ${\mathcal T}(H, \nu )$.
\end{rem}

3.2. Step 1. Let $ [\sigma _1  ] , [ \sigma _2 ] \in H_2(M, \Z )$ be homology
classes,  $\Sigma _j$ has genus $g_j$ ($j = 1, 2 $); 
 we can assume that both $\Sigma _1 $, $\Sigma _2 $  are 
 hyperbolic (otherwise the intersection pairing vanishes).
 Let $\psi _j: \pi _1(\Sigma _j) \to G$ be the
representation induced by $\sigma _j  $.

Now fix $j$ and set $\Sigma := \Sigma _j$ until the step 7 .

\bigskip
3.3. Step 2. 
The group $\psi ( \pi _1(\Sigma )) $ contains at least one loxodromic
element; hence we can construct  a pleated map 
$f^0 : \Sigma \to  M$ inducing $\psi : \pi _1(\Sigma ) \to G$ (see \cite{T}).

The pleated locus $\mathcal L$ of $f^0$ is a  geodesic lamination on $\Sigma $ .

Pick a maximal union $L_0$ of simple closed disjoint geodesics $\gamma $
on $\Sigma $ such that 
$$
0 < l_{\Sigma } (\gamma ) \le \mu \ .$$

\bigskip

3.4. Step 3. For every component $P_j \subset \Sigma  - L_0$ and $\mu > 0 $ 
 introduce the set 
$$
W_\mu (P_j) = \{ z \in P_j : InjRad(z) \le \mu /2 \}$$
 
Each ideal boundary component $\alpha \subset  \partial P_j $ has
orientation induced from $ P_j$ so we shall distinguish curves
$\alpha \subset L_0$ with different orientations but equal underlying sets. Set:

$W_\mu (\alpha , P_j) =   \{ z \in P_j :$ there exists a loop $
\beta _z $ on $P_j $ which is homotopic to $\alpha $ and passes
 through $z$, so
 that $l_{\Sigma }(\beta _z) \le \mu \} $. 

Then 
$$
W_\mu (P_j) = \bigcup _{ \alpha \subset \partial P_j } W_\mu (\alpha , P_j) $$

The Margulis constant $\mu $ is $\le $  the Margulis constant for 
$\H ^2$; therefore, for different boundary
 components we  obtain: $W_\mu (\alpha , P_j) \cap W_\mu 
(\beta ,
 P_j) = \emptyset $. Set  
$$ 
  \Sigma _\mu  = \{ z \in  \Sigma  : InjRad(z) \ge \mu /2 \}.$$
 Let $L_1 = L_0
- cl (\Sigma _\mu )$ and $P^0_j = P_j - W_\mu (P)$ for every $j$. Define
$$
\diam_0 (\Sigma _\mu ) = \sum_j \diam (P^0_j ) .$$

\begin{lem}\label{l8}
%LEMMA 8. 
$$
 \diam_0 (\Sigma _\mu ) \le (2g - 2)/\mu + \hbox{length~of} \partial \Sigma_\mu. $$
\end{lem}
\proof (Cf. \cite{Ab, Bo, T}.) We shall denote by $ l(\partial \Sigma _\mu )$
 the
 length of $\partial \Sigma_\mu $. Cover $\Sigma_\mu $ by a maximal set of
 disjoint discs $D(x , 2\mu )$. Then the number of these discs is
  $n < Area ( \Sigma ) /(4 \mu ^2)$
and every point $z \in \Sigma_\mu $ has the property:
$$
d(z, \partial (\Sigma_\mu  - \bigcup _{i=1}^n D_i (z, 2\mu )) \ < 2\mu$$

Hence, for every $z, w \in P_j^0 $ we have:
$d (z, w) \le 4n\mu + l(\partial \Sigma_\mu ) \le (2g - 2)/\mu + l(\partial
 \Sigma_\mu )$.
\qed

%REMARK 9. 
If $\alpha ^* \subset \partial \Sigma_\mu $ 
 then $ l(\alpha ^*)  \le 2\sinh
\mu $ by  Proposition \ref{p6}. Therefore we obtain:

\begin{cor}\label{c8}
%COROLLARY 8. 
In Lemma \ref{l8} we have:
$$
  \diam_0 \Sigma _\mu \le (2g-2)/\mu  + 6(g-1) \sinh \mu =  C_2(\mu , g ).  $$ 
\end{cor}

3.5. Step 4. 

Let $T_\gamma (\mu ) \subset M_{(0, \mu ]}$ 
be the Margulis tube whose fundamental group contains the
$G$-conjugacy class of
$\psi (\gamma )$, $\gamma \subset  L_1 $.

\begin{rem}\label{r10}
%REMARK 10. 
In general, $\pi _1 T_\gamma (\mu ) \ne \langle \psi (\gamma )\rangle$.
\end{rem}

 For every such geodesic $\gamma$ we have
 two (probably equal) components $P_i , P_j \subset \Sigma -
L_1$ adjacent to 
$\gamma$. Then $f^0 (W_\mu (\gamma , P_j) ) \subset T_\gamma (\mu ) $ \ 
$(k= i , j)$.

Choose points $x_k = x_{\gamma , k} \in \partial W_\mu (\gamma , P_k) 
\ (k= i , j) $ such that:

$f^0(x_i) , f^0(x_j) \in \Pi _t $ for some fiber of the canonical foliation of
$T_\gamma (\mu )$ associated with $\gamma $.

Let $\nu = \min \{ InjRad (f_0(x_k)), $ over all points $x_{\gamma , k}$
and all $\gamma \subset  L_1 \} $.

\begin{prop}\label{p77}
%PROPOSITION 7. 
$\nu = C_3(\mu , g)$ for some function $C_3$ which does not 
depend on the manifold $M$.
\end{prop}
\proof For every $i$ 
we have a point $o_i \in P_i - f^{-1}(M_{(0, \mu]})$ since $ \psi (\pi _1 (
 P_i ))$
is not virtually abelian.
Then $d_M(o_i , x_i) \le C_2(\mu , g)$.
Therefore (by Lemma \ref{l3}) 
$$
InjRad_M(f \ x_i)  \ge C_1(C_2(\mu , g), \mu ) = C_3 (\mu , g)$$
 \qed

%REMARK 11. 
\begin{rem}\label{r11}
Above we used the assumption that $\psi $ is a monomorphism.
\end{rem}

Let $\alpha ^*_k , \beta^*_k ,..., \omega^*_k $
 be the boundary components of $P^0_k$.
Then we can ``triangulate" $P^0_k$ so that: vertices of this ``triangulation"
$\Omega _f$
are $x_{\alpha } , x_{\beta },..., x_{\omega }$; lengths of edges of the
triangulation are bounded from above by $(6g-6)^2C_2(\mu , g)$ (see Figure \ref{F8} 
for a triangulation of the pair of pants). The triangles from this
triangulation will be  called ``short".

\bigskip

3.6. Step 5. Now, for each $k$ we map the triangulated surface $P_k^0$ to a
 p-g surface in $M$ by the new map 
$f : P^0_k \to M$ so that: for every edge $e$ of the triangulation we have: $f(e) \sim f^0(e)
 (\mathrm{rel} \ \partial e ) $.

Hence, $l_M(f(e)) \le l_{\Sigma }(e) \le (16g- 16)^2C_2(\mu , g)$. 

Now consider the thin part of $\Sigma $.

\bigskip

3.7. Step 6. Fix  $x_{\gamma , i} \ ,  \ x_{\gamma , j} $ lying on the
components 
$P_i \ , \ P_j$ 
adjacent to $\gamma \subset L_1 $.  Then connect their images 
$f(x_{\gamma , i}) \ ,  \ f( x_{\gamma , j} ) $ by the p-g annulus
$F_{\psi (\gamma )f(x_{\gamma , i})f( x_{\gamma , j} )}$ (see Definition \ref{d2}).
The boundary of $F= F_{\psi (\gamma )f(x_{\gamma , i})f( x_{\gamma , j} )}$
is equal to $f\gamma ^*_ i \cup f\gamma ^*_ j$. The annulus $F$ consists
of two geodesic triangles. These triangles will be called ``long" triangles
 ``sitting
 in" $T_\gamma (\mu )$. The  annulus $F$ itself will be called
``long p-g annulus sitting in $T_\gamma (\mu )"$. (See Remark \ref{r8}.)

Thus, we extended our map from $\Sigma _\mu $ to the p-g map
 $f: \Sigma \to M $ which is  
 homotopic to $\sigma $.

\begin{lem}\label{l9}
%LEMMA 9. 
Suppose that $d(p(A_{\psi \gamma }) , z) \le R + 2$ for some
$z \in f(\Sigma _\mu)$. Then 
$$
\diam (f \Delta ) \le 4 \sinh (\mu /2) /C_1(3R + 2 , \mu)$$
for every $f( \Delta ) $ sitting in $ T_\gamma (\mu ) $.
\end{lem}
\proof 
We have $d(z, f \ o_i) \le R , f ( o_i) \in M_{(\mu , \infty ]} \cap f(P^0),
d(f(o_i) , p(A_{\psi \gamma })) \le 3R + 2$ and so $l(\gamma ) \ge C_1(3R + 2, \mu )$
by Lemma \ref{l3}. This implies that
$$
d(p(A_{\psi \gamma }) , z) \le 2 \sinh (\mu
/2)/C_1(3R + 2 , \mu )$$
for every $z \in  T_{\gamma }(\mu ) $. Therefore for the  triangle 
$f( \Delta )$ sitting in $ T_\gamma (\mu )$ we have:
$$
\diam (\Delta ) \le 4 \sinh (\mu /2) /C_1(3R + 2 , \mu)$$
\qed

Our construction of the map $f$ is sufficiently flexible and given any two 
classes $\sigma _1 , \sigma _2 \in H_2(M , \Z )$ we can find transversal
 maps $f_i : \Sigma _i \to M$ representing these classes.

\bigskip

3.8. Step 7. The maps $f_i$ above will be called ``nice." Below is a summary of 
 properties of "nice" maps $f_i$: 

(1) $f_i$ are p-g with respect to some triangulations $\Omega _{f_i}$ of the
surfaces $\Sigma _i$. 
The number of triangles in $\Omega _{f_i}$ is $\le 16(g_i-1)$.

(2) In the triangulation $\Omega_{f_i} $ there are ``short" and ``long" triangles
(with respect to p-l metric induced by $f_i$). Namely, internal diameter
of ``short" triangles is 
$$
\le 2C_2(\mu , g_i) = R$$
and their union is $\Sigma
_{\mu i }$ so
that $\Sigma _i - \Sigma
_{\mu i} = W_{\mu i} $ is the union of pairwise disjoint nonhomotopic tubes.
 All vertices of $\Omega_{f_i} $ are contained in $\partial \Sigma _{\mu i}$.

(3) $f(0$-skeleton of $\Omega _{f_i} ) \subset M_{[\nu , \infty )}$ for
some $\nu = \nu (g_i)$.

(4) Every short triangle $\Delta \subset  P^0_{ki} \subset \Sigma _{\mu _i }$
 has bounded 
index with respect to any $\gamma ^* \subset \partial P^0_{ki} $.
 More precisely,

for every $z \in \dot \Delta $ and $\gamma ^* \subset \partial  \Sigma _{\mu i} $ passing
through $z$ we have: $l(\gamma ^* ) \le R $ (this is just the corollary
of (2)).

(5) Denote by $\Omega ^1 _{f_i}$  the 1-skeleton of $\Omega _{f_i}$.
 Every component $Q_{ji} \subset W _{\mu i}$  consists
of two ``long" triangles such that: if 
for some 
$h \in \pi _1(M)$ we have
 $$
 1 \ne  h^{-1} \psi _1( \pi _1 ( Q_{j1}) )h \cap \psi _2
( \pi _1 ( Q_{j2})) = \langle \gamma _j \rangle \subset \pi _1(M) $$
then $ \langle \gamma _j \rangle$ is a maximal  virtually abelian subgroup of  $\pi _1(M) $ and
 the maps 
$$
f_i : Q_{ji} - (\Omega ^1 _{f_i}) \to M$$
 can be lifted to the fundamental domain $\Phi _j \subset \H ^4$
of the group $\pi _1(T_{\gamma _j} (\mu )) $. The fundamental domain
$\Phi _j$ is bounded by pair of fibers of the canonical fibration of
$\H ^4$ corresponding to $ \langle \gamma _j \rangle$.

This property is the corollary of the condition MAX and Step 6.

(6) For every ``long triangle" $\Delta $ we have: $ \diam(f(\Delta ) 
\cap M_{[\nu , \infty )} ) \le C_4$ ; where 
$$
C_4 =  4R + 6 + 1/k(R, \nu ) + 4 \sinh (\mu /2 ) /C_1(3R +2 , 
\mu ) \}$$
this follows from
Lemma \ref{l5}  and Lemma \ref{l9}.

(7) Suppose that $d(p(A_{\psi \gamma }) , z) \le R + 2$ for some
$z \in f_i(\Sigma _{i \mu})$. Then 
$$
\diam (f_j \Delta ) \le 4 \sinh (\mu /2) /C_1(3R + 2 , \mu)$$
for every long triangle $f_j( \Delta ) $ sitting in $ T_\gamma (\mu ) $, 
for both $j= 1, 2$ (see  Lemma \ref{l9}).

\bigskip

3.9. Step 8. Now we can count 
the number of intersections
 $\#(f_1 (\Sigma _1) \cap f_2 (\Sigma _2))$.

(i) Consider intersections of "short" triangles. Pick a pair of such
 triangles $f_1 (\Delta _1 ) \subset f_1 (\Sigma _1)  , f_2(\Delta _2)
 \subset f_2
(\Sigma _2) $\ , let $\Delta ' _j \subset \H ^4$
 be the geodesic triangles covering them $(j= 1, 2)$. Then we are to estimate the
 number of 
$h \in G$ such that $h \Delta ' _1 \cap \Delta '_2 $ is non-empty. Recall that
$ \diam\Delta '_j \le R/2$ and both   $\Delta '_j $ contain points $o_j$
such that $Ir_G(o_j) \ge \mu /2 $. Therefore we can apply Lemma 2 to obtain:

$\#(f_1 (\Delta _1) \cap f_2 (\Delta _2)) \le 8 \exp ^3 (2R + \mu
 /2))/{\mu ^3}$. 

(ii) First consider the  case when $\Delta _1 $ is short while $\Delta _2 $ is long.

Suppose that $ h(z) \in h \Delta ' _1 \cap \Delta '_2  $. Then we can apply
Lemma \ref{l5} and the Property 7 of nice maps to obtain 
$$
d(z, \dot \Delta_2') \le \max \{ 4R + 6 + 1/( C_1 (R + 2 , \nu )), \ C_4 \}
 = C_5 \ .
$$
 Let $\{w_1 , w_2 , w_3 \} = \dot \Delta_2'$. 

Hence we obtain estimate in the same manner as in the case (i):
$$
\#(f_1 (\Delta _1) \cap f_2 (\Delta _2)) \le \# \{ h \in G : h( B(w_i , 2C_5))
 \cap \Delta _1' \ne \emptyset , i= 1, 2, 3 \}$$
$$ 
 \le 3 \exp ^3 (2C_5 + \nu /2))/
{\nu ^3}$$
since $ \diam(\Delta _1') \le 2R < C_5$.

(iii) Assume now that both $\Delta _1 , \Delta _2 $ are long. Denote by 
$ T_{\gamma _j}(\mu ) \subset M_{(0, \mu ]}$ those Margulis tubes
where $\Delta _j$ are sitting. 

(a) First count the
number of intersections that occur in $ T_\gamma (\mu )$ if 
 $\Delta _1 , \Delta _2 $ are sitting in the same 
$ T_\gamma (\mu ) \subset M_{(0, \mu ]}$. The influence of the fundamental
 group is
trivial (property 5) and here we have not more than $1$ intersections between
the  ``long"
triangles $f_1(\Delta _1 ), f_2( \Delta _2) $.

(b) If $ T_{\gamma _j}(\mu )$ are different then
$$
f_1(\Delta _1 ) \cap T_{\gamma _1}(\mu ) \cap
f_2( \Delta _2)  \cap T_{\gamma _2}(\mu ) = \emptyset $$ 

So consider intersections $(f_1(\Delta _1 ) - T_{\gamma _1}(\mu )) \cap
 f_2( \Delta _2) $ outside $T_{\gamma _2}(\mu )$ . However 
$\Delta _1  - f_1^{-1}
(T_{\gamma _1}(\mu )) $ is the union of 2 subsets each having diameter 
$ \le C_4 < C_5$ and therefore the number of  intersections is not more than
$$
6 \exp ^3 (2C_5 + \nu /2)/
{\nu ^3}$$
analogously to the case (ii). 

Hence the total estimate
in the case (iii) is
$$
6(12 (g-1))^2 \exp ^3 (2C_5 + \nu /2)/
{\nu ^3}$$ 
 where $g = \max \{g_1 , g_2 \}$. 

Direct calculations now show that the number 
$|< [\sigma _1 ], [\sigma _2 ] >|$ can be estimated as:

$$
|< [\sigma _1 ], [\sigma _2 ] >| \le 300 (g-1)^2 \exp (4000(g-1)/ \mu ) / \mu ^2$$

This concludes the proof of  Theorem \ref{t1}  under the condition MAX. \ \ \ \  \qed

\section{Proof of  Theorem \ref{t1}  in general}

\bigskip
4.1. The proof proceeds in the same way as in the section 3 until the step 6.

Step 6'. Fix  $x_{\gamma , i} \ ,  \ x_{\gamma , j} $ lying on the
components 
$P_i \ , \ P_j$ 
adjacent to $\gamma \subset L_1 $. Lift 
$x_{\gamma , i} \ ,  \ x_{\gamma , j} $ to points ${\tilde fx}_{\gamma , i} \ ,
  \ {\tilde fx}_{\gamma , j} $ in $ {\mathcal K}(H ,\mu )$ where  $ H$ is
the maximal  virtually abelian subgroup of $G$ which contains $<\psi (\gamma )>$.
Denote by $u' $ the projection $\varphi _{ H,\nu}(u)$ of the point
 $u \in \H ^4$ to the $
\partial {\mathcal K}(H ,\nu )$. For $k \in \{i, j\}$ we have:

(i) $d({\tilde fx}_{\gamma , k}, \psi (\gamma ){\tilde fx}_{\gamma , k}) \le
\mu$;

(ii) $Ir_G( {\tilde fx}_{\gamma , k}) \ge \nu /2$.

Therefore, according to Proposition \ref{p7}, either

(a) 
$d(\varphi {\tilde  fx}_{\gamma , k}, {\tilde fx}_{\gamma , k} )
\le C_+(\mu , \nu )$
or

(b) $\cosh d(\varphi {\tilde fx}_{\gamma , k}, A_{\psi (\gamma )}) \le
 \frac{2\sinh \mu /2}{  C_{-}(\mu , \nu )} $.

On another hand,
$ f (  x_{\gamma , k}) \in M_{(\nu , \infty ]}$.
 Therefore, if the case (b) holds then (by Lemma 3) we obtain the
lower estimates  
$$
C_1(\mathrm{arccosh}\left(\frac{2\sinh \mu /2}{ C_{-}(\mu , \nu )}\right), \nu) \le l(\psi \gamma ) \le l _\Sigma (\gamma )$$ 
and moreover:
$$
C_1(\mathrm{arccosh}\left(\frac{2\sinh \mu /2}{C_{-}(\mu , \nu )}\right), \nu) \le Ir_G(u)$$
for every $u \in A_{\psi (\gamma )}$.

Then we can consider the triangles which constitute the p-g annulus 
$$
F= F_{\psi (\gamma )f(x_{\gamma , i})f( x_{\gamma , j} )}$$
 as ``short triangles" and exclude $T_\gamma(\nu )$ from the consideration of the
thin part of $M$. Thus, let us assume that the alternative (b) does not hold.
Denote by $x''_{\gamma, k}$ the projections
 $p(\varphi {\tilde f}x_{\gamma, i})$

Then we construct the extended ruled film ${\tilde Q}= {\tilde Q}_{\gamma
(\varphi{\tilde fx}_{\gamma, i}) (\varphi{\tilde fx}_{\gamma, j}) } $; and connect
the loop 
$f^0(\gamma^*_ i )$ with $  \partial_{ \varphi {\tilde f}x_{\gamma, i}}
 p{\tilde Q}$
and 
$ f^0(\gamma^*_ j)$ with $\partial_{\varphi {\tilde f} x_{\gamma , j}} p{\tilde Q}$
by the p-g annuli 
$$
F_{\varphi {\tilde f} x_{\gamma , i},  f_0x_{\gamma , i}}\ , 
F_{\varphi {\tilde f} x_{\gamma , j},  f^0x_{\gamma , j}}\ $$
 (see Figure \ref{F9}). Now, instead of the
 long p-g 
annulus $F$ sitting in $T_\gamma (\mu ) $ (as in 3.7) take the union:

\begin{figure}[htbp] %  figure placement: here, top, bottom, or page
   \centering
   \includegraphics[width=4in]{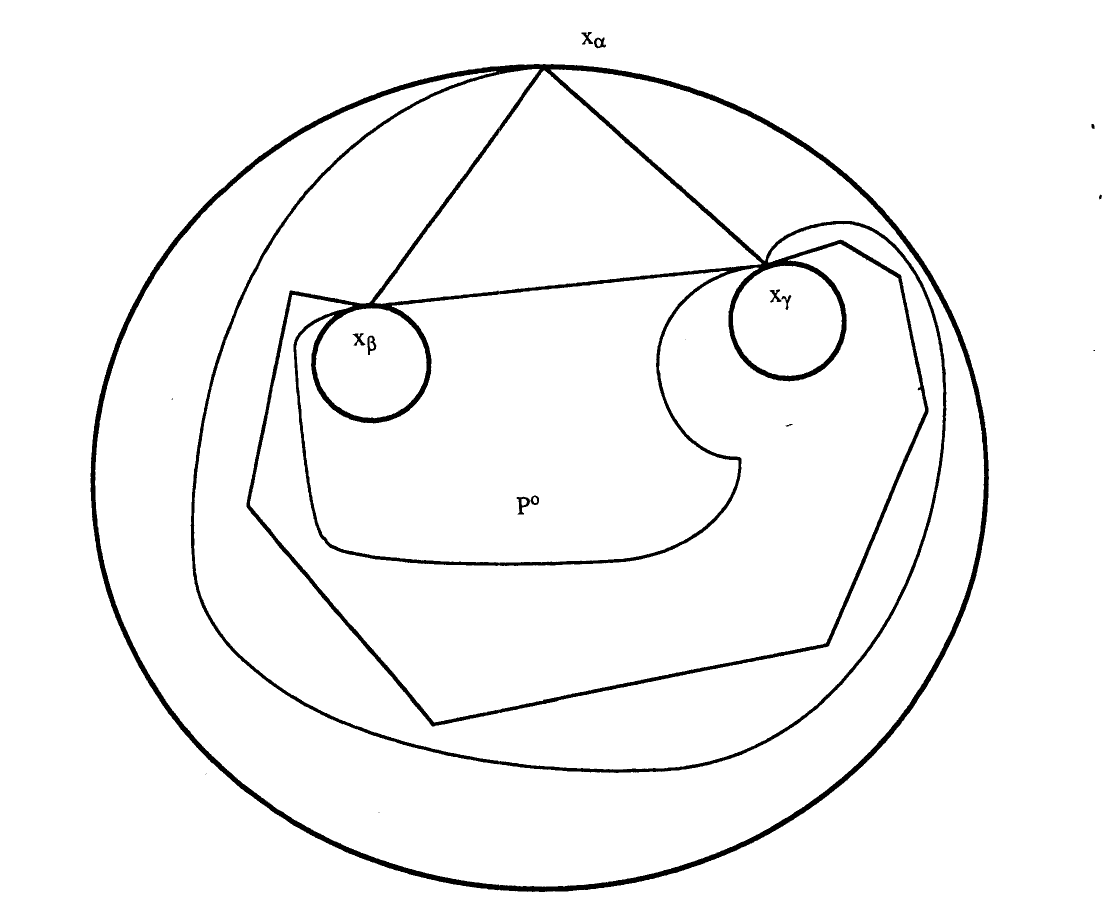} 
   \caption{}
   %\caption{example caption}
   \label{F9}
\end{figure}

$$
 p( {\tilde Q}_{\gamma \varphi 
{\tilde fx}_{\gamma , i} \varphi {\tilde fx}_{\gamma , j} })) \cup 
F_{ x''_{\gamma , i},  f^0x_{\gamma , i}} \cup 
F_{ x''_{\gamma , j},  f^0x_{\gamma , j}} .$$

In this way we extended 
the p-g map from $\Sigma _\mu$ to the new map $\phi $ defined on $\Sigma $.

4.2. Step 7'. The maps $\phi _i$ above will be called ``improved nice" maps.
 Restrictions of $\phi _i$ to $\Sigma _{i, \mu }$ have all properties of  nice
 maps. However
instead of ``long p-g annuli sitting in $T_\gamma (\mu ) $" we are  using
 extended ruled annuli (see Step 6').

4.3. Now we can estimate the intersection number 
$<\phi _1 (\Sigma _1) , \phi _2 (\Sigma _2)>$ in the same way as it was done in 3.9.

(i) Number of intersections between short triangles is estimated exactly as in
3.9(i).

(ii) Every extended ruled film $\tilde Q_k$ is contained in $2R$-neighborhood of 
the geodesic segment 
$$
  [{\tilde \phi }_k(x_{\gamma , i}) \ , 
 \ {\tilde \phi }_k(x_{\gamma , j})]   $$
($k= 1, 2$). 
Number of points of transversal intersection between any extended ruled film 
with geodesic plane is not more than $8$. Therefore again we can use Lemma \ref{l5} 
to obtain
$$
\# (p {\tilde Q _1} \cap f_2\Delta  ) \le  24 \exp ^3 (2C_5 + \nu /2))/
{\nu ^3}$$ 
for every short geodesic triangle  $f_i\Delta $ (see 3.9(ii)).

\begin{rem}\label{r12}
%REMARK 12. 
If the loxodromic alternative holds for the segment   $[{\tilde \phi }_1(x_{\gamma , i}) \ , 
 \ {\tilde \phi }_1(x_{\gamma , j})]$ holds  then 
the loxodromic alternative holds for $[{\tilde \phi }_2(x_{\gamma , i}) \ , 
 \ {\tilde \phi }_2(x_{\gamma , j})]$ too.
\end{rem}

(iii) The algebraic intersection number between two extended ruled annuli
 that belong to one and the same Margulis tube was estimated in Theorems
 \ref{t4}, \ref{t5}. It is not more than 
$$
\max \{ N(R, \nu ) , N'(R, \nu ) , R/ \nu \}$$
This finishes the proof of  Theorem \ref{t1}  in the general case. \ \ \ \  \qed

\section{Appendix}

\begin{lem}\label{l11}
%LEMMA 11. 
Let $p_1$,  $p_2$ be  closed geodesics on hyperbolic surface
$F$. Then:
$$
\#(p_1 \cap \ p_2 ) \le \exp(l_1 + l_2 + 1)$$
where $\#(. \cap .)$ is the number of points of intersection, $l_i$ is
the length of $p_i$.
\end{lem}

\proof Denote by  $G \subset Isom(\H ^2)$ the fundamental group
of $F$.   Denote by $g_i$ representative
of $p_i$ in the deck-transformation group $G$; and let $q_i$ be the
axis of $g_i$ . Let $f_i$ be a segment of
 the length $l_i$ on  $q_i$ . Without loss of
generality we can assume that there is $x \in f_1 \cap f_2 $ such that
 $Ir_G(x) \ge 1$
(since 2 can be used as  the Margulis
constant in $\H ^2$). Then
$$
<p_1, \ p_2> \le \# \{ h \in G : h (B(x, l_1 + l_2)) \cap B(x, l_1 + l_2)
\ne \emptyset \}$$
Now we can apply 2-dimensional version of Lemma 2 to obtain that
$$
 \# \{ h \in G : h (B(x, l_1 + l_2)) \cap B(x, l_1 + l_2) \ne \emptyset
 \} \le \exp (l_1 + l_2 + 1)$$

Lemma is follows. \qed 

%\newpage

%\centerline {\bf References}


\begin{thebibliography}{BLP05}

\bibitem{Ab} W. Abikoff, {\em The Real Analytic Theory of  Theichm\"uller Space},
Lecture Notes in Mathematics, vol 820. Springer Verlag, 1980.

\bibitem{A} M. Anderson, {\em Metrics of negative curvature
on vector bundles}, Proceedings of AMS, 99 (1987), p. 357--363. 


\bibitem{Bo} F. Bonahon, {\em Boutes des varietes hyperboliques dee dimension 3}, Ann.
of Math. 124 (1986) p. 71--158.

\bibitem{B} B.H. Bowditch, {\em Geometrical finiteness for hyperbolic groups}, Journal of Functional
Anal., 113 (1993), p. 245--317.

\bibitem{CT} J. Carlson, D. Toledo, {\em Harmonic mappings of K\"ahler manifolds to
 locally symmetric spaces}, Publications of IHES, vol. 69 (1989) p. 173-- 201.

\bibitem{F} W. Fulton,  ``Intersections Theory'', Series of Modern Surveys in
Mathematics, Springer, 1984. 

\bibitem{GT} M. Gromov, W. Thurston, {\em Pinching constants for hyperbolic
manifolds}, Inv. Math. v. 89, F. 1 (1987) p. 1-- 12.

\bibitem{GLT} M. Gromov, H.B. Lawson, W. Thurston, {\em Hyperbolic structures on
4-manifolds and flat conformal structures on 3-manifolds},
Math. Publications of IHES, vol. 67 (1988) p. 27--45.
  
\bibitem{HI} E. Heintze, H. Im Hof,  {\em Geometry of horospheres}, J. Diff. Geom.,
1977, 12, p. 481--491.

\bibitem{Ki} R. Kirby, {\em Problems in low-dimensional topology}. In: Algebraic and
 Geometric Topology. (Proc. Symp. Pure Math. v. 22) AMS, Providence, 1978, 
p. 273-- 312

\bibitem{Kr} P. Kronheimer, {\em Embedded Surfaces in 4-manifolds}, Proceedings of ICM-90,
vol. 1, p. 529-- 538.
 
\bibitem{Ku1} N. Kuiper, {\em Hyperbolic 4-manifolds and tessellations (variations
on (GLT))}, Math. Publications  of IHES, vol. 68 (1988), p. 47-- 76

\bibitem{Ku2} N. Kuiper, {\em Fairly symmetric hyperbolic manifolds }, In:
``Geometry and Topology of Submanifolds, II'' (1990), p. 165- 203.
World Sci. Publisher, Singapore, New Jersey, London, Hongkong

\bibitem{Ku3} N. Kuiper, 1988 (unpublished)

\bibitem{Ka1} M. Kapovich, { \em Flat conformal structures on 3-manifolds. The existence
problem, I.} Siberian Math. Journal, vol 30 (1989)  N 5 p. 60--73.


\bibitem{Ka2} M. Kapovich, {\em Flat conformal structures on 3-manifolds (survey)},
 In: Proceedings of International Conference dedicated to A.Maltsev,
Novosibirsk 1989, Contemporary Mathematics, 1992, Vol. 131.1, p. 551--570. 

\bibitem{Ka3} M. Kapovich, { \em Deformation spaces of  flat conformal structures},
 In: ``Proceedings of Soviet--Japanese Topology Symposium held in
 Khabarovsk 1989.''  Answers and Questions in General Topology, vol. 8 (1990) p. 253--264.

\bibitem{Ka4} M. Kapovich, { \em Flat conformal structures on 3-manifolds. I} J. Diff.
Geom. Vol. 38, N 1, (1993) p. 191--215.

\bibitem{GKL} W. Goldman, M. Kapovich and B. Leeb, {\em Complex hyperbolic surfaces homotopy-equivalent to a Riemann surface}, Comm. in Analysis and Geom., Vol. 9, N 1 (2001), p. 61--96.

\bibitem{Mo} N. Mok, {\em K\"ahler Geometry of Arithmetic Varieties}, In: ``Several Complex
Variables and Complex Geometry'' 
 (Proc. Symp. Pure Math. v. 52 , part 2) AMS, Providence, 1989, p. 335-- 396

\bibitem{Mor} J. Morgan, {\em Group actions on trees and the compactification of the
 space of classes of} $SO(n, 1)$ {\em representations}, Topology, 1986,
 V. 25, N 1, p. 1--33.

\bibitem{MS} G. Mostow, Y.-T. Siu, {\em A compact K\"ahler surface of negative curvature
not covered by the ball}, Ann. Math. 112 (1980), p. 321-- 360.

\bibitem{SY} Y-T. Siu, S.-T. Yau, {\em K\"ahler manifolds of finite volume},  Ann.
Math. Stud., N 102, 1982, p. 363-379

\bibitem{T} W. Thurston,  ``Geometry and Topology of 3-manifolds'', Lecture Notes,
Princeton University, 1978.

\bibitem{To1} D. Toledo,  {\em Harmonic maps from surfaces to
		certain K\"ahler manifolds,} 
		Math. Scand. {\bf 45} (1979), p. 13-- 26

\bibitem{To2} D. Toledo, 
                {\em Representations of surface groups on
		complex hyperbolic space,}
                J. Diff. Geo. {\bf 29} (1989), p. 125--133

\end{thebibliography}
\end{document}